%% file: Abelian_Varieties_With_Homotheties.tex
\documentclass[10pt, reqno, intlimits]{amsart}
\usepackage{amsthm,amsmath,amssymb,amsopn,amsfonts,amscd,amsbsy,amsrefs}
\usepackage{mathtools}
\usepackage{url}
\usepackage{ifthen}
\usepackage[all, cmtip]{xy}
\usepackage{pdfpages}
\usepackage{graphics}
\usepackage{mathrsfs}
\usepackage{hyperref}
\hypersetup{
	colorlinks   = true,
	urlcolor     = blue,
	linkcolor    = blue,
	citecolor   = red
}

\usepackage{tikz-cd}

\title{Abelian varieties with homotheties}

\author{F\'elix Baril Boudreau}
\address{Departamento de Matem\'aticas\\
CUNEF Universidad, Madrid, Espa\~na}
\email{felix.baril@cunef.edu}

\author{Sebastian Petersen}
\address{Fachbereich Mathematik und Naturwissenschaften\\
Universität Kassel, Kassel, Germany}
\email{petersen@mathematik.uni-kassel.de}



\date{\today}
\input{math}
\newcommand{\ab}{\mathrm{ab}}
\newcommand{\conn}{\mathrm{conn}}
\newcommand{\pr}{\mathrm{pr}}
\newcommand{\Aut}{\mathrm{Aut}}
\newcommand{\sep}{\mathrm{sep}}
\newcommand{\Zar}{\mathrm{Zar}}
\newcommand{\Fr}{\mathrm{Fr}}
\newcommand{\etabar}{\overline{\eta}}

\newcommand{\HA}[1]{\HHH_{#1}(\Zz_\ell)}
\newcommand{\HAl}{\HA{A,\ell}}
\newcommand{\HAL}{\prod_{\ell\in L} \HAl}
\newcommand{\HAlq}{\HHH_{A,\ell}^{(q)}}
\newcommand{\HALq}{\HHH_{A, L}^{(q)}}
\newcommand{\HALqj}{\HHH_{A, L}^{(q^j)}}
\newcommand{\CAl}{C_{A,\ell}}
\newcommand{\CAL}{C_{A,L}}
\newcommand{\orho}{\overline{\rho}}

\newcommand{\condH}{\hyperlink{conditionH}{(H)}}
\newcommand{\defcondH}{\hypertarget{conditionH}{(H)}}

\begin{document}
		\maketitle
		
		\begin{abstract}
			In 1980, Bogomolov proved that the $\ell$-adic Galois image of an Abelian variety over a number field contains an open subgroup of the homotheties. In 1986, Serre refined that result, showing that the homotheties coming from the Galois side were uniformly bounded in $\ell$ inside the homotheties of the codomain of the representation. In a 2007 paper, Zarhin shows that Bogomolov's conclusion can fail in non-trivial cases for global function fields. We prove a suitable analogue of Serre's results over finitely generated fields of positive characteristic. As a consequence, we prove a new case of the 1978 conjecture by Geyer and Jarden.
		\end{abstract}

		\section{Introduction}
		
		\subsection{Background and previous work}
		
		Let $K$ be a field of characteristic $\chara(K)$ with separable closure $K^\sep$ in a fixed algebraic closure $\oK$ and let $\Gal(K)$ be its absolute Galois group. Let $A/K$ be an Abelian variety of dimension $g_A=\dim(A) \geq 1$. Let $\Ll$ be the set of prime numbers different from $\chara(K)$. For each $\ell \in \Ll$ we have the $\Zz_\ell$-adic and $\Qq_\ell$-adic {\bf Tate modules}
		$$
		T_\ell(A) := \varprojlim_{n \geq 1} A[\ell^n](\oK) \ \mbox{and}\ V_\ell(A) := T_\ell(A) \otimes_{\Zz_\ell} \Qq_\ell
		$$
		of $A$. The group schemes $A[\ell^n]$ are finite \'etale over $K$. Hence $A[\ell^n](\oK)=A[\ell^n](K^\sep)$ and we thus get a continuous action of 
		$\Gal(K)$ on $T_\ell(A)$ yielding an $\ell$-adic Galois representation 
		$$\rho_{A,\ell}: \Gal(K)\to \GL_{T_\ell(A)}(\Zz_\ell),$$
		where $\GL_{T_\ell(A)}$ is the $\Zz_\ell$-group scheme such that $\GL_{T_\ell(A)}(R)=\Aut_R(T_\ell(A)\otimes_{\Zz_\ell} R)$ for all commutative $\Zz_\ell$-algebras $R$. For a non-empty subset $L\subset \Ll$ we let 
		$$\rho_{A, L}: \Gal(K)\to \prod_{\ell\in L} \GL_{T_\ell(A)}(\Zz_\ell)$$ be the group morphism induced by the $\rho_{A,\ell}$'s with $\ell \in L$.
		For every separable extension $K'/K$ write $G_{A,L}(K')$ for the image $\rho_{A,L}(\Gal(K'))$. 
		A general principle in our notation is that $G_{A,\ell}(K')=G_{A, \{\ell\}}(K')$, $\rho_{A, \ell}=\rho_{A, \{\ell\}}$ and so forth.
		By the $\ell$-adic variant of Cartan's theorem, $G_{A,\ell}(K)$ is a compact subgroup of 
		$\GL_{V_\ell(A)}(\Qq_\ell)$, hence an $\ell$-adic Lie group. We denote by $\mfg_{A, \ell}$ the Lie algebra of $G_{A,\ell}(K)$ and by $\GGG_{A,\ell}$ the Zariski closure of $G_{A,\ell}(K)$ in $\GL_{V_\ell(A)}$. 
		We denote by $\HHH_{A,\ell}$ the group of homotheties in  $\GL_{V_\ell(A)}$, that is, a copy of $\Gg_m$ embedded diagonally into $\GL_{V_\ell(A)}$. Then $\HHH_{A,\ell}$ is the center of the algebraic group $\GL_{V_\ell(A)}$. We conveniently define
		$$\GGG_{A,\ell}(\Zz_\ell):=\GGG_{A,\ell}(\Qq_\ell)\cap \GL_{T_\ell(A)}(\Zz_\ell)\ \mbox{and}\ \HAl:=\HHH_{A,\ell}(\Qq_\ell)\cap \GL_{T_\ell(A)}(\Zz_\ell).$$
		Finally we let $\CAL(K'):=G_{A, L}(K')\cap \HAL$. We now turn to two natural questions. First, how much of the center $\HAl$ of $\GL_{T_\ell(A)}(\Zz_\ell)$ is contained in the $\ell$-adic Galois image $G_{A,\ell}(K)$? Second, how far is $G_{A, \ell}(K)$ from the $\Qq_\ell$-points of its Zariski closure $\GGG_{A,\ell}$? Seminal work of Bogomolov and Serre answers the former question in the number field case as follows.
		
		\begin{thm}{\bf (Bogomolov \cite{bogomolov_1980}*{Th\'eor\`eme~2~and~Corollaire~1}, \cite{bogomolov_1980_2}*{Theorem~3})}\label{Bogomolov_Theorem}
			Let $A$ be an Abelian variety of dimension $g_A\ge 1$ over a number field $K$.
			\begin{enumerate}
				\item[(a)] For each $\ell \in \Ll$, the group $G_{A,\ell}(K)$ is open in $\GGG_{A,\ell}(\Qq_\ell)$.
				\item[(b)] For each $\ell \in \Ll$, the Lie algebra $\mfg_{A,\ell}$ is algebraic.
				\item[(c) ] For each $\ell \in \Ll$, the Lie algebra $\mfg_{A,\ell}$ contains $\Qq_\ell \cdot \Id_{V_\ell(A)}$. 
			\end{enumerate}
		\end{thm}
		Part (a) implies that $[\GGG_{A, \ell}(\Zz_\ell): G_{A, \ell}(K)]<\infty$ and Part (c) implies that we have
		$[\HAl:\CAl(K)]<\infty$. 
		Serre subsequently quantified Theorem \ref{Bogomolov_Theorem}:
		
		\begin{thm}{\bf (Serre \cite{Serre_1986}*{p.~92,~Th\'eor\`eme 2}, \cite{Serre_1986}*{p.~93,~Corollaire})}\label{Serre_Theorem_A}
			Let $A$ be an Abelian variety of dimension $g_A\ge 1$ over a number field $K$. 
			\begin{enumerate}
				\item[(a)] The function $f_A: \Ll\to \Nn^*, \ell\mapsto [\HAl : \CAl(K)]$ is bounded. 
				\item[(b)] There is an integer $c(A/K) \geq 1$ depending on $A/K$, such that $G_{A,\Ll}(K)$ contains all $c(A/K)$-th powers of elements of  
				$\prod_{\ell\in\Ll} \HAl$.
			\end{enumerate}
		\end{thm}
		
		These theorems about the existence of homotheties inside the image of Galois have since then seen a wide range of applications, for example in the context of the Manin-Mumford conjecture \cite{Oesterle_1985}, in work on the Geyer-Jarden conjecture \cite{JJ}, in monodromy computations such as \cite{2003-Banaszak_Gajda_Krason}, 
		in Galois cohomology \cite{1971-Serre_Congruence_II},  \cite{2005-Larsen}, or, indirectly, in work on detecting linear dependence via reduction maps
		\cite{2005-Banaszak_Gajda_Krason} or in relation with the support problem \cite{larsen-support}. Zywina addressed the natural question of how $G_{A,\ell}(K)$ sits in its Zariski closure.
		
		\begin{thm}{\bf (Zywina \cite{zywina-open}*{Theorem 1.2 (b)})}\label{zywina-thm}
			Let $A$ be an Abelian variety over a number field $K$. 
			There is a finite set $\Sigma=\Sigma(A/K)$ depending on $A/K$ such that the function 
			$$\Ll\setminus \Sigma\to \Nn^*,\ \ell\mapsto [\GGG_{A, \ell}(\Zz_\ell): G_{A, \ell}(K)]$$ is bounded by a constant that depends only 
			on the dimension $g_A$ of $A$. 
		\end{thm}
		
		From now on, $K$ is a {\em finitely generated field of positive characteristic $p$} and $A/K$ an Abelian variety. It is natural to ask how much can be recovered from the above theorems in that case. 
		
		\begin{defi}
			We say that $A$ satisfies condition \defcondH{}  if 
			for each $\ell \in \Ll$, the  Lie algebra $\mfg_{A,\ell}$ contains $\Qq_\ell \cdot \Id_{V_\ell(A)} $.
		\end{defi} 
		
		Thus, the Abelian varieties satisfying condition \condH~are those that satisfy the conclusion of Part (c) of Bogomolov's Theorem \ref{Bogomolov_Theorem}. 
		Zarhin points out in \cite{zar-hom} that an ordinary elliptic curve over a finite field does not satisfy condition \condH. 
		Moreover Zarhin exhibits absolutely simple Abelian varieties over a global function field $K$ for which condition \condH~fails and not isogenous over $\oK$ to an Abelian variety over a finite field \cite{zar-hom}*{Theorem~1.4,~Example~4.3}.
		
		In Proposition  \ref{Int:Equiv} below it will become transparent that the conclusions of Theorem \ref{Bogomolov_Theorem},  
		of Theorem \ref{Serre_Theorem_A} and of Theorem \ref{zywina-thm} fail drastically for those Abelian varieties A/K that do not satisfy condition \condH. On the flip side, the following important theorem of Zarhin guarantees the existence of many Abelian varieties satisfying condition \condH.
		
		\begin{thm} {\bf (Zarhin, \cite{zar77}*{Corollary~1}, \cite{zar-finchar14}*{Section~6})} If the center of the $\Qq$-algebra $\End_{\ol{K}}(A_{\ol{K}}) \otimes_\Zz \Qq$ is a product of totally real fields, then $A$ satisfies \condH.
		\end{thm}
		
		The main purpose of our paper is to establish {\em suitable} analogues of Theorem \ref{Bogomolov_Theorem},  
		of Theorem \ref{Serre_Theorem_A} and of Theorem \ref{zywina-thm}
		for those Abelian varieties over a finitely generated field of positive characteristic that satisfy condition \condH.  The {\em direct} analogues turn out to fail in positive characteristic, even in the presence of condition \condH, as we shall see in
		Proposition \ref{int:failure}. 
		
		\subsection{Main results} 
		
		Let $K$ be a finitely generated extension of a finite field $k$ and $p=\chara(k)$. 
		Let $A/K$ be an Abelian variety of dimension $g_A\ge 1$. We continue to denote by $\Ll$ the set of all prime numbers different from $p$. 
		Our first results offer several equivalent characterizations of condition \condH~and prove that these conditions do not depend on $\ell$. 
		Let $\ZZZ_{A,\ell}$ be the center of the connected component of the identity $\GGG_{A,\ell}^\circ$ of $\GGG_{A,\ell}$. 
		
		\begin{intprop} (Lemmas \ref{lemm:Hl} and \ref{lemm:H2})\label{Int:Equiv} Let $\ell\in\Ll$. The following conditions are equivalent. 
			\begin{enumerate}
				\item[(a)] The Lie algebra $\mfg_{A,\ell}$ contains $\Qq_\ell \cdot \Id_{V_\ell(A)}$. 
				\item[(b)] The Lie algebra $\mfg_{A,\ell}$ is algebraic. 
				\item[(c)] The identity component $\ZZZ_{A,\ell}^\circ$ of $\ZZZ_{A,\ell}$ equals $\HHH_{A,\ell}$.
				\item[(d)] The algebraic group $\ZZZ_{A, \ell}$ has dimension $1$.
				\item[(e)] The group $\CAl(K)$ is open in the compact group $\HAl$. 
				\item[(f)] The group $\CAl(K)$ is infinite. 
				\item[(g)] The group $G_{A, \ell}(K)$ is open in the compact group $\GGG_{A,\ell}(\Zz_\ell)$. 
			\end{enumerate}
			Moreover, if these equivalent conditions hold for one $\ell\in \Ll$, then they hold for all $\ell\in\Ll$ and $A$ satisfies condition \condH. 
		\end{intprop}
		
		If $A$ does not satisfy condition  \condH, then Proposition \ref{Int:Equiv} tells us that $[\GGG_{A, \ell}(\Zz_\ell):G_{A, \ell}(K)]=\infty$ and 
		$[\HAl: \CAl(K)]=\infty$ for all $\ell\in \Ll$. 
		However, even if we assume condition \condH, there is a further obstruction to a direct analogue of Theorems \ref{Serre_Theorem_A}
		and \ref{zywina-thm} as we now explain.
		
		Existence and basic properties of the Weil pairing force $\GGG_{A, \ell}$ to be contained in the group $\GSp_{V_\ell(A)}$ of symplectic similitudes with respect to the Weil pairing (see Section \ref{Section:weil_and_cyclotomic}). Let $q$ be the cardinality of the exact constant field of $K$. For any (always non-empty) subset $L$ of $\Ll$, let $q_L = (q_\ell)_{\ell \in L}$ be the image of $q$ under the unique ring morphism $\Zz\to \prod_{\ell\in L}\Zz_\ell$. For each integer $j \geq 1$, note that $q_L^j = (q_\ell^j)_{\ell \in L} \in \prod_{\ell\in L}\Zz_\ell^\times$ because $p\notin \Ll$ and $L\subset \Ll$. We then let $\langle q_L^j \rangle$ be the closure of the subgroup $(q_L^j)^\Zz$ of $\prod_{\ell\in L}\Zz_\ell^\times$ generated by $q_L^j$. For each $\ell \in L$, let $m_{A,\ell} : \GSp_{V_\ell(A)} \to  \Gg_{m,\Qq_\ell}$ be the algebraic group morphism that sends a symplectic similitude to its multiplier. Writing likewise $m_{A,\ell}$ for the resulting group morphism on its $\Qq_\ell$-points, let $m_{A,L} := \prod_{\ell \in L} m_{A,\ell} : \prod_{\ell \in L} \GSp_{V_\ell(A)}(\Qq_\ell) \to \prod_{\ell \in L} \mathbb{G}_{m,\Qq_\ell}(\Qq_\ell)$ for these later $m_{A,\ell}$'s. Finally we put for all positive integers $j$
		\begin{align*}
			&\HALqj:=\left\{h\in \HAL: m_{A,L}(h)\in \langle q_L^j \rangle\right\}\ \mbox{and}\\
			&\GGG_{A, L}^{(q^j)} := \left\{g\in \prod_{\ell\in L}\GGG_{A, \ell}(\Zz_\ell): m_{A,L}(g)\in \langle q_L^j \rangle\right\}. 
		\end{align*}
		
		\begin{intprop} (Corollary \ref{rema:topo_gen}, Proposition \ref{Prop:Unbounded_Homothety_Function} and Corollary \ref{coro:notbounded2}) \label{int:failure}
			\begin{enumerate}
				\item[(a)] For every subset $L$ of $\Ll$ there are inclusions  $G_{A, L}(K)\subset \GGG_{A, L}^{(q)}$ and $\CAL(K)\subset \HALq$. 
				\item[(b)] The function $\Ll\to \Nn^*, \ell\mapsto [\HHH_{A, \ell}(\Zz_\ell): \HHH_{A, \ell}^{(q)} ]$ is unbounded.
								\item[(c)] The function $\Ll\to \Nn^*, \ell\mapsto [\GGG_{A, \ell}(\Zz_\ell): \GGG_{A, \ell}^{(q)}]$ is unbounded. 
			\end{enumerate}
		\end{intprop}
		
		Proposition \ref{int:failure} has the following immediate Corollary. 
		
		\begin{intcoro} (Corollary \ref{coro:notbounded})
			\begin{enumerate}
				\item[(a)] The function $\Ll\to \Nn^* \cup \{\infty\}:\ \ell \mapsto [\HAl:\CAl(K)]$  is unbounded.
				\item[(b)] The function $\Ll\to \Nn^* \cup \{\infty\}:\ \ell \mapsto [\GGG_{A,\ell}(\Zz_\ell): G_{A,\ell}(K)]$  
				is unbounded 
			\end{enumerate}
		\end{intcoro}
		
		One is thus naturally led to ask whether the indices  $[\HAlq : \CAl(K)]$ and 
		$[\GGG_{A, \ell}^{(q)}: G_{A, \ell}(K)]$ can be bounded independently of  $\ell$ in case $A/K$ satisfies condition \condH. 
		
		\begin{intthm} (Corollary \ref{coro:withconstants3})
			Suppose $A/K$ satisfies condition \condH. There is an integer $c_1=c_1(g_A) \geq 1$ depending only on $g_A$ and a
			finite subset $\Sigma\subset \Ll$ such that the following hold. 
			\begin{enumerate}
				\item[(a)]  We have $[\GGG_{A,\ell}^{(q)}:G_{A, \ell}(K)]\le c_1$ for all $\ell\in\Ll\setminus \Sigma$.	
				\item[(b)]  We have $[\HHH_{A,\ell}^{(q)}:C_{A, \ell}(K)]\le  c_1$ for all $\ell\in\Ll\setminus \Sigma$.
				\item[(c)] We have $g^{c_1!}\in G_{A, \Ll\setminus \Sigma}(K)$ for all $g\in 
				\GGG^{(q)}_{A, \Ll\setminus \Sigma}$. 
				\item[(d)] We have $h^{c_1!} \in C_{A, \Ll\setminus \Sigma}(K)$ for all $h\in 
				\HHH^{(q)}_{A, \Ll\setminus \Sigma}$. 
			\end{enumerate}
		\end{intthm}
		
		Up to replacing $c_1$ by a constant that depends on more than just $g_A$, we can get rid of $\Sigma$.
		
		\begin{intcoro} (Corollary \ref{coro:withconstants2})
			Suppose $A/K$ satisfies condition \condH.  
			\begin{enumerate}
				\item[(a)] The function $\Ll \to \Nn^*, \ell\mapsto [\GGG_{A,\ell}^{(q)}:G_{A, \ell}(K)]$ is bounded.
				\item[(b)] The function $\Ll \to \Nn^*, \ell\mapsto [\HHH_{A,\ell}^{(q)}:C_{A, \ell}(K)]$ is bounded. 
			\end{enumerate}
		\end{intcoro}
		
		\begin{intcoro}  (Corollary \ref{coro:withconstants2}) Suppose $A/K$ satisfies condition \condH. 
			\begin{enumerate}
				\item[(a)] There is an integer $c_2 \geq 1$ such that $g^{c_2}\in G_{A, \Ll}(K)$ for all $g\in 
				\GGG^{(q)}_{A, \Ll}$. 
				\item[(b)] There is an integer $c_3 \geq 1$ such that $g^{c_3}\in C_{A, \Ll}(K)$ for all $g\in 
				\HHH^{(q)}_{A, \Ll}$. 
			\end{enumerate}
		\end{intcoro}
		
		Let us briefly discuss the main ingredients in the proofs of these results. By Zarhin's work, the $\Gal(K)$-module $V_\ell(A)$ is semisimple (see \cite{moret-bailly_1985}*{p.~244,~Corollaire~2.4} and \cite{zar-finchar75}*{Corollary~2.5}) and the center $\mfz_{A,\ell}$ of 
		the Lie algebra $\mfg_{A,\ell}$ is a $1$-dimensional $\Qq_\ell$-vector space (see \cite{zar77}*{Theorem 2}). 
		
		A further important ingredient is the observation that - regardless of condition \condH~- the Zariski closure $\GGG_{A,\ell}$ contains the algebraic group 
		$\HHH_{A,\ell}$ of homotheties. We establish this in Lemma \ref{lemm:CenterHomotheties} with the help of information on the eigenvalues of the Frobenii coming from the Weil conjectures.  
		
		The equivalent statements of Proposition \ref{Int:Equiv} are then established using the above in combination with the theory of $\ell$-adic Lie groups, of their Lie algebras and of their Zariski closures. Moreover, the $\ell$-independence in Proposition \ref{Int:Equiv} uses an important theorem of Chin (see \cite{chin_2004}*{Theorem 1.4 (ii), p.~726}).
		The basic properties of the Weil pairing and the behaviour of the cyclotomic character of $K$ yield Proposition \ref{int:failure}~(a), Proposition \ref{int:failure}~(b) is shown by a direct argument that we owe to Arno Fehm. To derive Proposition \ref{int:failure}~(c) from Proposition \ref{int:failure}~(b) one again relies on the information that 
		$\GGG_{A,\ell}\supset \HHH_{A,\ell}$ regardless of condition \condH. 
		
		Let $\GGG_{A,\ell}^\geo$ be the Zariski closure of $G_{A,\ell}(\overline{k}K)$ in $\GL_{V_\ell(A)}$ and $\KKK_{A,\ell}$ the kernel of the restriction of the algebraic group morphism $m_{A,\ell}$ to $\GGG_{A,\ell}^\circ$. Under condition \condH~the group scheme $\KKK_{A,\ell}/\GGG_{A,\ell}^{\geo, \circ}$ is shown to be finite \'etale of rank at most $2$. The final ingredient is a theorem of Cadoret, Hui and Tamagawa which implies that, up to replacing $K$ by a finite extension, the index $[\GGG_{A,\ell}^\geo(\Zz_\ell):\rho_{A,\ell}(\Gal(\ok K))]$ is bounded independently of $\ell$.  
		
		\subsection{An application}
		Let $K$ be a finitely generated field and $e\in\Nn^*$. Equip $\Gal(K)$ with the unique normalized Haar measure $\mu$ such that $\mu(\Gal(K))=1$. Then extend this measure to $\Gal(K)^e$ and denote it also $\mu$.
		Every $\eta\in \Gal(K)$ extends uniquely to an automorphism of the algebraic closure $\ol{K}$ of $K^\sep$ denoted again by $\eta$. For every extension $N$ of $K$ contained in $\ol{K}$ and every 
		$\sigma:=(\sigma_1,\cdots, \sigma_e)$ in $\Gal(K)^e$ let $N(\sigma)$ be the fixed field in $N$ of $\sigma_1,\cdots, \sigma_e$. In particular $\ol{K}(\sigma)$ is the 
		maximal purely inseparable extension of $K^\sep(\sigma)$. 
		In 1978 Geyer and Jarden stated the following conjecture and solved it for elliptic curves. 
		
		\begin{conj}[Geyer--Jarden \cite{GJ-ELL}]\label{Geyer_Jarden_Conjecture} Let $K$ be a finitely generated field of characteristic $p\ge 0$. For all $\sigma\in \Gal(K)^e$ outside a subset of $\Gal(K)^e$ of measure zero, each Abelian variety $A$ of dimension $g_A\ge 1$ defined over $\oK(\sigma)$ has the following properties.
			\begin{enumerate}
				\item[(a)] If $e= 1$, then there are infinitely many prime numbers $\ell$ with $A(\ol{K}(\sigma))[\ell]\neq 0$. 
				\item[(b)] If $e\ge 2$, then there are only finitely many  prime numbers $\ell$ with $A(\ol{K}(\sigma))[\ell]\neq 0$. 
				\item[(c)] If $e\ge 1$, then the $\ell$-primary part $A(\ol{K}(\sigma))[\ell^\infty]$ is finite for all prime numbers $\ell$. 
			\end{enumerate}
		\end{conj}
		
		Part (a) in case $\chara(K)=0$ was proved in \cite{jp-gj} building on \cite{gj-torsion} and \cite{zywina}. Part (c) is meanwhile known to be true in full generality \cite{JJ}. The same paper proves Part (b) in the case $\chara(K)=0$. More precisely, in \cite{JJ}*{The homothety condition on a field K (2a)} $K$ is said to satisfy the \textit{homothety condition} if for every Abelian variety $A/K$ of positive dimension the index 
		$[\HAl: C_{A, \ell}(K)]$ is bounded independently from $\ell$. They show Part (b) holds when $K$ satisfies this condition (see \cite{JJ}*{Proposition 3.3}). But by 
		Theorem \ref{Serre_Theorem_A}, Proposition \ref{int:failure} and Corollary \ref{coro:notbounded} a finitely generated field $K$ satisfies the homothety condition if and only if $\chara(K)=0$. Thus Part (b) remains open for $p>0$. In \cite{2013-Arias-de-Reyna_Gajda_Petersen}*{Theorem 4.1} 
		Part (b) was proved for principally polarized Abelian varieties with big monodromy.
		
		Combining our main results with ingenious ideas from the paper \cite{JJ} of Jacobson and Jarden, we can establish the following contribution to the conjecture. 
		
		\begin{intthm}  (Theorem \ref{gjthm}) Let $K$ be a finitely generated field of positive characteristic, $A/K$ an Abelian variety and $e\ge 3$ an integer. Assume $A$ satisfies condition \condH. For all $\sigma\in \Gal(K)^e$ outside a subset of $\Gal(K)^e$ of measure zero there are only finitely many  prime numbers $\ell$ with $A(\oK(\sigma))[\ell]\neq 0$.
		\end{intthm}
		
		The case $e = 2$ is more subtle. It will be part of an ongoing work. 
		
		\subsection{Notation and Conventions}
		Let $k$ be a finite field of characteristic $\chara(k):= p$ and let $\Ll$ be the set of primes numbers but $p$. Let $K/k$ be a finitely generated field extension and $A/K$ an Abelian variety.
		
		Let $\ell \in \Ll$. For a compact subgroup $G$ of $\GL_{V_\ell(A)}(\Qq_\ell)$ we denote by $\Lie(G)$ its Lie algebra and by $G^\Zar$ the Zariski closure of $G$ in $\GL_{V_\ell(A)}$. For a closed subgroup scheme $\MMM$ of $\GL_{V_\ell(A)}$ we denote by $\Lie(\MMM)$ its Lie algebra, by $\MMM^\circ$ its connected component of identity and we put $\MMM(\Zz_\ell):=\MMM(\Qq_\ell) \cap \GL_{T_\ell(A)}(\Zz_\ell)$. Reductive and semisimple group schemes are assumed to have connected fibres over the base. 
		We define for $K'/K$ separable and (always non-empty) $L\subset \Ll$ the group $G_{A,L}^\geo(K')=\rho_{A,L}(\Gal(\ol{k}K'))$, and let $\mfg_{A,\ell}^\geo=\Lie(G_{A,\ell}^\geo(K))$ and $\GGG_{A,\ell}^\geo=G_{A,\ell}^\geo(K)^\Zar$. For a locally compact topological group $G$ and $X\subset G$ we denote by $\langle X \rangle$ the topological closure of the subgroup of $G$ generated by $X$. In our notation $\Nn$ and $\Nn^*$ are respectively the set of non-negative and positive integers.
		
\subsection*{Acknowledgements}
The project started when F.B.B. was a
postdoctoral fellow at the Universit\'e du Luxembourg. This research was
funded in part by the Luxembourg National Research Fund (FNR), grant
reference 501100001866-17921905. 
For the purpose of open access, and in
fulfilment of the obligations arising from the grant agreement, the
authors have applied a Creative Commons Attribution 4.0 International
(CC BY 4.0) licence to any Author Accepted Manuscript version arising
from this submission. While F.B.B. was a CICMA postdoctoral fellow at Concordia University 
from April to August 2026, he received funding from the NSERC grant DG-155635-2019 of Chantal David.
The project was continued and completed when
F.B.B. was Assistant Professor at CUNEF Universidad. S.P. thanks Gebhard Böckle and Wojciech Gajda for 
discussions and explanations concerning homotheties some years ago. S.P. acknowledges a visit at Adam Mickiewickz University Pozna\'n in autumn 2026, where part of this work was done. We both thank Arno Fehm for help with Lemma \ref{Lemma:Not_Bounded}.

		\section{Preliminaries on monodromy groups}\label{Section:weil_and_cyclotomic}
		Let $k$ be a finite field of characteristic $p$, $K/k$ a finitely generated extension and $A/K$ an Abelian variety of dimension $g_A\ge 1$. Let $q$ be the cardinality of the exact constant field $k'$ of $K$. For $\ell\in\Ll$ we define $\Zz_\ell(1):=\plim_{n\in \Nn^*} \Gg_m[\ell^n](\oK)$ and $\Qq_\ell(1):=\Zz_\ell(1)\otimes_{\Zz_\ell} \Qq_\ell$. Both are written multiplicatively in our notation and are free of rank $1$ over their respective rings $\Zz_\ell$ and $\Qq_\ell$.
		By \cite{milne-abvar}*{\S 16, pp. 131-132}, the Weil pairing $e_\ell: V_\ell(A)\times V_\ell(A^\vee)\to \Qq_\ell(1)$ is a non-degenerate $\Gal(K)$-equivariant pairing for every $\ell\in \Ll$. Choose an ample line bundle $\LLL\in \mathrm{Pic}(A)$ and let $\lambda=\varphi_\LLL: A\to A^\vee$ be the corresponding isogeny \cite{milne-abvar}*{p. 119} and in fact a polarization. It thus induces an isomorphism $\lambda: V_\ell(A)\to V_\ell(A^\vee)$. Following \cite{milne-abvar}*{\S 16, p. 132} one constructs for every $\ell\in \Ll$ a pairing
		$$
		e_{A,\ell} := e_\ell^\LLL: V_\ell(A)\times V_\ell(A)\to \Qq_\ell(1),\ (x,y)\mapsto e_\ell(x, \lambda(y))
		$$
		that is non-degenerate and $\Gal(K)$-equivariant. By \cite{milne-abvar}*{Lemma 16.2~(e)} it is also alternating. For each $\ell \in \Ll$, let $\varepsilon_\ell: \Gal(K) \to \Zz_\ell^\times$ be the $\ell$-adic cyclotomic character and let $\GSp_{V_\ell(A)}$ be the general symplectic $\Qq_\ell$-algebraic group with respect to the pairing $e_{A,\ell}$. For each non-empty subset $L \subset \Ll$, let also $\varepsilon_L := \prod_{\ell \in L} \varepsilon_\ell : \Gal(K) \to \prod_{\ell \in L} \Zz_\ell^\times$.
		
		\begin{lemm} \label{lemm:mdet} Let $L$ be a subset of $\Ll$. The following hold.
			\begin{enumerate}
				\item[(a)] For all $\ell\in\Ll$ we have $G_{A,\ell}(K)\subset \GSp_{V_\ell(A)}(\Qq_\ell)$ and $\GGG_{A,\ell}\subset \GSp_{V_\ell(A)}$. 
				\item[(b)] For all $\sigma\in\Gal(K)$ we have 
				$(m_{A, L}\circ \rho_{A, L})(\sigma)=\varepsilon_L(\sigma)$.
				\item[(c)] We have $\varepsilon_L(\Gal(K))=\langle q_L\rangle$. 
				\item[(d)] Let $(a_\ell)_{\ell \in L} \in \prod_{\ell\in L} \Gg_{m,\Qq_\ell}(\Qq_\ell)$ and let $(a_\ell \Id_{V_\ell(A)})_{\ell\in L}$ be the corresponding element of $\prod_{\ell\in L} \HHH_{A, \ell}(\Qq_\ell)$. Then $m_{A,L}((a_\ell \Id_{V_\ell(A)})_{\ell\in L})= (a_\ell^{2})_{\ell \in L}$.
			\end{enumerate}
		\end{lemm}
		
		\begin{proof}
			Let $\ell \in \Ll$. By $\Gal(K)$-equivariance of $e_{A,\ell}$, we have
			$$
			e_{A,\ell}(\rho_{A, \ell}(\sigma) x, \rho_{A, \ell}(\sigma) y)=e_{A, \ell}(x,y)^{\varepsilon_\ell(\sigma)}
			$$
			for all $\sigma \in \Gal(K)$ and $x, y \in V_\ell(A)$. Hence $G_{A,\ell}(K)\subset \GSp_{V_\ell(A),e_{A,\ell}}(\Qq_\ell)$ and Part (a) follows by taking the Zariski closures of these two sets. But $m_{A,\ell}$ is the multiplier map on $\GSp_{V_\ell(A),e_{A,\ell}}$. Therefore,
			$$
			e_{A,\ell}(\rho_{A, \ell}(\sigma) x, \rho_{A, \ell}(\sigma) y) = e_{A,\ell}(x,y)^{m_{A,\ell}(\rho_{A,\ell}(\sigma))}
			$$
			for all $\sigma \in \Gal(K)$ and $x, y \in V_\ell(A)$. As $e_{A,\ell}$ is non-degenerate, $m_{A,\ell}(\rho_{A,\ell}(\sigma))$ equals $\varepsilon_\ell(\sigma)$ for all $\sigma \in \Gal(K)$. This implies Part (b). If $k'$ is the exact constant field of $K$, then the restriction map $r: \Gal(K) \to \Gal(k')$ is surjective. Moreover the group morphism $\varepsilon_L: \Gal(K)\to \prod_{\ell\in L} \Zz_\ell^\times$ factors through $r$. Indeed, if $\varepsilon_{k', L}: \Gal(k')\to \prod_{\ell\in L} \Zz_\ell^\times$ is the cyclotomic character of the finite field $k'$, then $\varepsilon_L=\varepsilon_{k', L}\circ r$. Now, let $q$ be the cardinality of $k'$. The profinite group $\Gal(k')$ is topologically generated by the Frobenius automorphism $F: \ol{k'}\to \ol{k'},\ F(x)=x^q$. Hence $\varepsilon_L(\Gal(K))=\varepsilon_{k',L}(\Gal(k'))$ is topologically generated by 
			$\varepsilon_{k', L}(F)$, and $\varepsilon_{k', L}(F)=q_L$ and Part (c) follows.
			
			Finally, if $a :=(a_\ell)_{\ell\in L} \in \prod_{\ell\in L} \Qq_\ell^\times$, for each $\ell \in L$, if $x_\ell,y_\ell\in V_\ell(A)$,
			$$
			e_{A,\ell}(a_\ell x_\ell,a_\ell y_\ell)
			=e_{A, \ell}(x_\ell,y_\ell)^{a_\ell^2}
			$$
			by the linearity of the pairing $e_{A, \ell}$. 
			Thus $m_{A,\ell}(a_\ell\Id_{V_\ell(A)})=a_\ell^2$
			for every $\ell\in L$, and this implies
			$m_{A,L}\bigl((a_\ell\Id_{V_\ell(A)})_{\ell\in L}\bigr)=(a_\ell^2)_{\ell\in L}$, proving Part (d).
		\end{proof}
		
		\begin{coro}\label{rema:topo_gen} 
			We have $G_{A, L}(K)\subset \GGG_{A, L}^{(q)}$ and $C_{A, L}(K)\subset \HHH_{A, L}^{(q)}$. 
		\end{coro}
		
		\begin{proof} We have $G_{A, \ell}(K)\subset \GGG_{A, \ell}(\Zz_\ell)=\GGG_{A, \ell}(\Qq_\ell)\cap \GL_{T_\ell(A)}(\Zz_\ell)$ for all $\ell\in\Ll$. By Lemma \ref{lemm:mdet}, $m_{A,L}(G_{A,L}(K)) = \varepsilon_L(\Gal(K)) = \langle q_L \rangle$. Thus, $G_{A,L}(K) \subset \GGG_{A, L}^{(q)}$. So, $C_{A, L}(K)=G_{A,L}(K)\cap \HAL \subset \GGG_{A, L}^{(q)} \cap \HAL \subset \HHH_{A, L}^{(q)}$.
		\end{proof}
		
		\begin{fact}\label{remark:semisimple} The group $G_{A,\ell}(K)$ acts semisimply on $V_\ell(A)$ by \cite{moret-bailly_1985}*{p.~244,~Corollaire~2.4} and \cite{zar-finchar75}*{Corollary~2.5}, and it follows from \cite{serr-dris}*{Proposition~2~(i),~(ii)} that $\GGG_{A,\ell}^\circ$ is reductive and from Clifford's Theorem \cite{curtis}*{p.~343,~Theorem~49.2} that every closed normal subgroup of $G_{A,\ell}(K)$ acts semisimply on $V_\ell(A)$. In particular 
			$G_{A,\ell}^\geo(K)$ acts semisimply on $V_\ell(A)$. 
		\end{fact}
		
		\begin{fact} \label{remark:derived} The algebraic group $\GGG_{A,\ell}^{\geo, \circ}$ is semisimple  and is the derived group of the reductive group $\GGG_{A, \ell}^\circ$ by \cite{del-weil2}*{Corollaire~1.3.9~and~Th\'eor\`eme~3.4.1~(iii)} and \cite{bgp-var}*{Corollary~1.3}, respectively.
		\end{fact}
		
		\begin{lemm} There is a smooth connected $k$-variety $S$ with function field $K$. 
		\end{lemm}
		\begin{proof} The field $k$ is perfect. Let $t=(t_1,\cdots, t_r)$ be a separating transcendence base for $K/k$. The affine space $\Aa^r_k$ has function field $k(t)$. Let $S$ be the normalization of $\Aa_k^r$ in $K$. The generic fibre of the finite surjective morphism $f: S\to \Aa^r_k$ is $\Spec(K)$. As $K/k(t)$ is separable, it follows that the generic fibre of $f$ is \'etale over $k(t)$. Hence there is a non-empty open subscheme $U$ of $\Aa^r_k$ such that the restriction $f_U: S_U\to U$ is \'etale. It follows that $S_U$ is a smooth connected $k$-variety with function field $K$ and we conclude after replacing $S$ by $S_U$. 
		\end{proof}
		
		\begin{fact} \label{fact:compatible}
			Let $S$ be a smooth connected $k$-variety with function field $K$ and geometric generic point $\overline{\eta} \to S$. 
			After replacing $S$ by one of its non-empty open subschemes,  the Abelian variety $A$ extends to an Abelian scheme $\pi: \AAA \to S$ \cite{milne-abvar}*{Remark~20.9}. 
			For $\ell\in \Ll$, the representation $\rho_{A,\ell}: \Gal(K)\to \GL_{V_{\ell}(A)}(\Qq_\ell)$ 
			is unramified at all codimension-one points of $S$ \cite{serre-tate_1968}*{Theorem~1}. By purity of the branch locus \cite{grothendieck_1968}*{Expos\'e X, Th\'eor\`eme~3.4}, $\rho_{A,\ell}$ 
			factors through the \'etale fundamental group $\pi_1(S,\overline{\eta})$. For every closed point $s\in S$, let $k(s)$ be its residue field and $\Phi_s\in\Gal(k(s))$ be the arithmetic Frobenius, characterized by $\Phi_s(x)=x^{\#k(s)}$ for all $x\in\ol{k(s)}$. 
			The image $\Fr_{A, s,\ell}$ of the Frobenius element $\Phi_s\in \Gal(k(s))$ under the 
			composite group morphism $\Gal(k(s))\to \pi_1(S,\overline{\eta}) \xrightarrow{\rho_{A,\ell}} \GL_{V_{\ell}(A)}(\Qq_\ell)$ is well-defined up to conjugation. Its characteristic polynomial has coefficients in $\Zz$, is independent of $\ell\in\Ll$ and the eigenvalues of $\Fr_{A, s,\ell}$ are pure of weight $1$ \cite{milne-abvar}*{Theorem~19.1}. 
		\end{fact}
		
		\begin{rema} \label{rema:chin_applies} The collection $(\rho_{A,\ell}: \pi_1(S,\etabar) \to \GL_{V_{\ell}(A)}(\Qq_\ell))_{\ell \in \Ll}$ of representations thus corresponds to a semisimple pure $\Qq$-compatible system $(\LL_\ell)_{\ell\in\Ll}$ of lisse $\Qq_\ell$-sheaves on $S$ in the sense of \cite{chin_2004}*{2.4}.
		\end{rema}
		
		\begin{fact} \label{fact:fingen} For any $\ell \in \Ll$ and integer $n \geq 1$, the profinite group $\GL_n(\Zz_\ell)$ has an open subgroup whose rank, in the sense of \cite{dixon2003analytic}*{Proposition~3.11,~Definition 3.12}, is $n^2$ by \cite{dixon2003analytic}*{Theorem~5.2}. So every closed subgroup of $\GL_n(\Zz_\ell)$ is topologically finitely generated \cite{dixon2003analytic}*{p.~89}.
		\end{fact}
		
		\begin{lemm}\label{Group_Theoretic_Lemma} Let $n\geq 1$ be an integer, $L$ a finite subset of $\Ll$ and $G$ a closed subgroup of $\prod_{\ell\in L} \GL_n(\Zz_\ell)$. Then the Frattini subgroup $\Phi(G)$ is open in $G$. 
		\end{lemm}
		
		\begin{proof} For each $\ell \in L$, let $\pi_\ell: \GL_n(\Zz_{\ell})\to \GL_n(\Ff_{\ell})$ be the reduction modulo $\ell$ map and 
			$\pr_{\ell}: \prod_{\ell\in L} \GL_n(\Zz_{\ell})\to \GL_n(\Zz_\ell)$ the projection on the $\ell$-th factor. The kernel $P_{\ell}$ of $\pi_\ell$ is open in 
			$\GL_n(\Zz_\ell)$ and pro-$\ell$. It follows that $N:=G \cap \prod_{\ell \in L} P_\ell$ is open in $G$. Each $\pr_\ell(N)$ is a subgroup of $P_\ell$, hence pro-$\ell$. Thus, for distinct $\ell_1, \ell_2 \in L$, $\pr_{\ell_1}(N)$ and $\pr_{\ell_2}(N)$ share no common finite quotient. By \cite{serre-ind}*{Lemme 2}, we have $N = \prod_{\ell \in L} \pr_\ell(N)$. Each $\pr_\ell(N)$ is a subgroup of $\GL_n(\Zz_\ell)$, hence topologically finitely generated by Fact \ref{fact:fingen}. Therefore, $\Phi(G)$ is open in $G$ \cite{Serre_1997}*{p.~148,~Proposition~3.6~(iv)~implies~(ii)}.
		\end{proof}
		
		\begin{lemm}Let $\ell\in\Ll$. 
			\begin{enumerate}
				\item[(a)] The canonical group morphism $G_{A,\ell}(K)\to \GGG_{A,\ell}/\GGG_{A,\ell}^\circ(\Qq_\ell)$ is surjective, the canonical group morphism
				$\GGG_{A,\ell}(\Qq_\ell)/\GGG_{A,\ell}^\circ(\Qq_\ell) \to \GGG_{A, \ell}/\GGG_{A,\ell}^\circ(\Qq_\ell)$ is bijective and the finite \'etale $\Qq_\ell$-group scheme $\GGG_{A,\ell}/\GGG_{A,\ell}^\circ$
				is constant. 
				\item[(b)]The canonical group morphism $G_{A,\ell}^\geo(K)\to \GGG_{A,\ell}^\geo/\GGG_{A,\ell}^{\geo, \circ}(\Qq_\ell)$ is surjective, the canonical group morphism
				$\GGG_{A,\ell}^\geo(\Qq_\ell)/\GGG_{A,\ell}^{\geo, \circ}(\Qq_\ell)\to \GGG_{A, \ell}^\geo/\GGG_{A,\ell}^{\geo, \circ}(\Qq_\ell)$ is bijective and the finite \'etale $\Qq_\ell$-group $		\GGG_{A,\ell}^\geo/\GGG_{A,\ell}^{\geo, \circ}$ is constant. 
			\end{enumerate}
		\end{lemm}
		
		\begin{proof} Let $\GGG=\GGG_{A,\ell}$ (respectively $\GGG=\GGG_{A,\ell}^\geo)$ and $G=G_{A,\ell}(K)$ (respectively $G=G_{A,\ell}^\geo(K)$). 
			Then $\pi: \GGG\to\GGG/ \GGG^\circ$ is the projection and $i: \GGG(\Qq_\ell)/\GGG^\circ(\Qq_\ell)\to \GGG/\GGG^\circ(\Qq_\ell)$ the canonical group morphism. 
			The group morphism $i$ is injective and the group $G$ is Zariski dense in $\GGG$. Hence $\pi(G)$ is Zariski dense in 
			the finite \'etale group scheme $\GGG/\GGG^\circ$. It follows that $\GGG/\GGG^\circ$ is constant and that $\pi: \GGG(\Qq_\ell)\to \GGG(\Qq_\ell)/\GGG^\circ(\Qq_\ell)\buildrel i\over \longrightarrow\GGG/\GGG^\circ(\Qq_\ell)$ is surjective. Thus
			$i$ is surjective.
		\end{proof}
		The following result is a slight adaption of results by Larsen and Pink. 
		
		\begin{lemm}\label{lemm:kernelindep} 
			\begin{enumerate}
				\item[(a)] The group $\rho_{A,\ell}^{-1}(\GGG_{A,\ell}^\circ(\Qq_\ell))$ does not depend on $\ell\in\Ll$. 
				\item[(b)]The group $\rho_{A,\ell}^{-1}(\GGG_{A,\ell}^{\geo, \circ}(\Qq_\ell))\cap \Gal(\ol{k}K)$ does not depend on $\ell\in\Ll$. 
			\end{enumerate}
		\end{lemm}
		
		\begin{proof} For $\ell\in\Ll$ consider the compositions $\rho_{A,\ell}^\conn: \Gal(K)\buildrel \rho_{A,\ell}\over \longrightarrow \GGG_{A,\ell}(\Qq_\ell)\to \GGG_{A,\ell}/\GGG_{A, \ell}^\circ(\Qq_\ell)$
			and $\rho_{A,\ell}^{\geo, \conn}: \Gal(\ol{k}K)\buildrel \rho_{A,\ell}\over \longrightarrow \GGG^\geo_{A,\ell}(\Qq_\ell)\to \GGG^\geo_{A,\ell}/\GGG_{A, \ell}^{\geo, \circ}(\Qq_\ell)$. We have
			to prove that $\ker(\rho_{A,\ell}^\conn)$ and $\ker( \rho_{A,\ell}^{\geo, \conn})$ do not depend on $\ell\in\Ll$. 
			By \cite{lar-pi0}*{Examples 6.3 and 6.6, Proposition 6.14}, for $\ell, \ell'\in\Ll$, we have
			$\rho_{A,\ell}^{-1}(\GGG_{A,\ell}^\circ(\Qq_{\ell}))=\rho_{A,\ell'}^{-1}(\GGG_{A,\ell'}^\circ(\Qq_{\ell'})) $ and 
			the groups  $\im(\rho_{A,\ell}^\conn)=\GGG_{A, \ell}(\Qq_\ell)/\GGG_{A,\ell}^\circ(\Qq_\ell)$ and 
			$\im(\rho_{A,\ell'}^\conn)=\GGG_{A, \ell'}(\Qq_\ell)/\GGG_{A,\ell'}^\circ(\Qq_\ell)$ are canonically isomorphic. Hence, Part (a) is true.
			
			Part (b) is trivial if $\mathrm{trdeg}(K/k)=0$ and established in \cite{lar-pi}*{Proposition 2.2} if $\mathrm{trdeg}(K/k)=1$. We deduce Part (b) for $\mathrm{trdeg}(K/k)\ge 2$ by the following specialization argument. Up to replacing $k$ by the exact constant field of $K$, we can assume that $S/k$ is geometrically connected. 
			Let $\ell_0, \ell_1 \in\Ll$ with $\ell_0\neq \ell_1$. Consider the group morphism
			$$
			\rho: \pi_1(S_{\ol{k}})\xrightarrow{(\rho_{A,\ell_0},\rho_{A,\ell_1})} \GGG_{A,\ell_0}^\geo(\Qq_{\ell_0}) \times \GGG_{A,\ell_1}^\geo(\Qq_{\ell_1})
			$$ 
			induced by $\rho_{A,\ell_0}$ and $\rho_{A,\ell_1}$. 
			By Lemma \ref{Group_Theoretic_Lemma} the Frattini subgroup $\Phi(\rho(\pi_1(S_{\ol{k}})))$ is open in $\rho(\pi_1(S_{\ol{k}}))$. Thus, by 	
			\cite{2023_Gajda_Petersen}*{Corollary~4.2}, there is a 	
			smooth connected curve $C$ on $S_{\ol{k}}$ and an immersion $i: C\to S_{\ol{k}}$ such that 
			$\rho(i_*(\pi_1(C))=\rho(\pi_1(S_{\ol{k}}))$ (*). This implies $\rho_{A, \ell_j}(i_*(\pi_1(C)))=\rho_{A,\ell_j}(\pi_1(S_{\ol{k}}))$ and $\rho_{A, \ell_j}(i_*(\pi_1(C)))^\Zar=\GGG_{A, \ell_j}^\geo$
			for $j\in \{0,1\}$. We have continuous group morphisms
			$$
			\pi_1(C)\buildrel i_*\over \longrightarrow \pi_1(S_{\ol{k}}) 
			\buildrel \rho_{A, \ell_j}|\pi_1(S_{\ol{k}})\over \longrightarrow \GGG_{A,\ell_j}^\geo(\Qq_{\ell_j}) \buildrel q_{\ell_j} \over \longrightarrow
			\GGG_{A, \ell_j}^\geo(\Qq_{\ell_j})/\GGG_{A, \ell_j}^{\geo, \circ}(\Qq_{\ell_j}),
			$$
			where $q_{\ell_j}$ is the canonical projection.
			There is a finite extension $k_0/k$ such that the curve $C$ and the immersion $i$ descend to a curve $C_0$ and an immersion $i_0: C_0\to S_{k_0}$ over $k_0$ and moreover $\rho_{A, \ell}\circ i_{0,*}$ is isomorphic to the $\ell$-adic Galois representation $\rho_{A', \ell}$ attached to the Abelian scheme $A':=\AAA\times_S C_0$ over $C_0$, where $\AAA/S$ is the Abelian scheme from Fact \ref{fact:compatible}. By \cite{lar-pi}*{Proposition 2.2} there is an isomorphism $\varphi: \GGG_{A, \ell_0}^\geo/\GGG_{A, \ell_0}^{\geo, \circ}(\Qq_{\ell_0})\to \GGG_{A, \ell_1}^\geo/\GGG_{A, \ell_1}^{\geo, \circ}(\Qq_{\ell_1})$ such that $\varphi\circ q_{\ell_0}\circ (\rho_{A, \ell_0}|\pi_1(S_{\ol{k}}))\circ i_*=q_{\ell_1}\circ (\rho_{A, \ell_1}|\pi_1(S_{\ol{k}}))\circ i_*$, and by (*) it follows that $\varphi\circ q_{\ell_0}\circ (\rho_{A, \ell_0} \pi_1(S_{\ol{k}}))=q_{\ell_1}\circ (\rho_{A, \ell_1}|\pi_1(S_{\ol{k}}))$, which implies (b). 
		\end{proof}
		
		\begin{defi} We let $K_A^\conn$ be the fixed field of $\rho_{A,\ell}^{-1}(\GGG_{A,\ell}^\circ(\Qq_\ell))$ in $K^\sep$ for some (hence, by Lemma \ref{lemm:kernelindep}, for every) $\ell\in\Ll$. Let $k_A^\conn$ be the exact constant field of $K_A^\conn$. We let $K_A^{\geo, \conn}$ be the fixed field of  $\rho_{A,\ell}^{-1}(\GGG_{A,\ell}^{\geo, \circ}(\Qq_\ell))\cap \Gal(\ol{k}K)$ in $K^\sep$ for some, hence every, $\ell\in\Ll$. Finally let $K_A^{(3,5,7)}$ be the compositum of the fields $K(A[\ell])$ for $\ell\in \{3,5,7\}\setminus\{p\}$.
		\end{defi}
		
		\begin{lemm} \label{lemm:Kconnbound} We let $\gamma=\gamma(2g_A):= \#\GL_{2g_A}(\Ff_3)\cdot \#\GL_{2g_A}(\Ff_5)\cdot\#\GL_{2g_A}(\Ff_7)$. 
			\begin{enumerate}
				\item[(a)] We have $K_A^\conn\subset  K_A^{(3,5,7)}$.
				\item[(b)] We have $K_A^{\geo, \conn}\subset \ol{k} K_A^{(3,5,7)}$. 
				\item[(c)] The index $[K_A^\conn:K]$ divides $\gamma$ and the index $[K_A^{\geo, \conn}:\ol{k}K]$ divides $\gamma$. If $k'$ is the exact constant field of $K$, then
				$[k_A^\conn:k']$ divides $\gamma$.
			\end{enumerate}
		\end{lemm}
		
		\begin{proof} Let $\Lambda=\{3,5,7\}\setminus\{p\}$ and let $E=K_A^{(3,5,7)}$. Then $E$ is the fixed field of the kernel of composition
			$$\Gal(K)\buildrel( \rho_{A,\ell})_{\ell \in \Lambda}\over \longrightarrow \prod_{\ell\in\Lambda} \GL_{T_\ell(A)}(\Zz_\ell)\to  \prod_{\ell\in\Lambda} \GL_{T_\ell(A)}(\Ff_\ell)\cong
			\prod_{\ell\in\Lambda} \GL_{2g_A}(\Ff_\ell),
			$$
			$[E:K]= \# \Gal(E/K)$ divides $\gamma$ (*) and $\rho_{A, \ell}(\Gal(E))$ is a closed subgroup of the pro-$\ell$-group $\ker(  \GL_{T_\ell(A)}(\Zz_\ell)\to  \GL_{T_\ell(A)}(\Ff_\ell))$ for $\ell\in \Lambda$. The group $\rho_{A,\ell}^{-1}(\GGG_{A,\ell}^\circ(\Qq_\ell))\cap \Gal(E)$ respectively $\rho_{A,\ell}^{-1}(\GGG_{A,\ell}^{\geo, \circ}(\Qq_\ell))\cap \Gal(\ol{k}E)$) 
			does not depend on $\ell\in\Ll$ by Lemma \ref{lemm:kernelindep}. Let $I_\ell$ be the image of $\Gal(E)$ under the group morphism
			$$
			\Gal(E)\xrightarrow{\rho_{A,\ell}}
			\GGG_{A,\ell}(\Qq_\ell)\longrightarrow
			\GGG_{A,\ell}(\Qq_\ell)/\GGG_{A,\ell}^\circ(\Qq_\ell),
			$$
			and let $\oI_\ell$ be the image of $\Gal(\ol{k}E)$ of the composition
			$$
			\Gal(\ol{k}E)\xrightarrow{\rho_{A,\ell}}
			\GGG_{A,\ell}^\geo(\Qq_\ell)\longrightarrow
			\GGG_{A,\ell}^\geo(\Qq_\ell)/\GGG_{A,\ell}^{\geo,\circ}(\Qq_\ell).
			$$
			Because the groups $\rho_{A,\ell}^{-1}(\GGG_{A,\ell}^\circ(\Qq_\ell))\cap \Gal(E)$ respectively $\rho_{A,\ell}^{-1}(\GGG_{A,\ell}^{\geo, \circ}(\Qq_\ell))\cap \Gal(\ol{k}E)$) 
			do not depend on $\ell\in\Ll$, 
			the same holds for the cardinalities $\#I_\ell$ and $\#\oI_\ell$ of these images. Furthermore for all $\ell\in\Lambda$ these cardinalities are powers of $\ell$. As $\#\Lambda\ge 2$ this implies $\#I_\ell=1$ and $\#\oI_\ell=1$, from which Parts (a) and (b) follow. Part (c) is immediate from (*) and Parts (a) and (b).
		\end{proof}
		
		For a subalgebra $\mfg$ of $\Lie(\GL_{V_\ell(A)}(\Qq_\ell))$ we denote by $\mfg^{\alg}$ its algebraic hull (see \cite{DG}*{II,~\S~6,~no.~2,~Definition 2.4}). We recall that 
		$\mfg_{A,\ell}=\Lie(G_{A,\ell}(K))$ and $\mfg_{A,\ell}^\geo=\Lie(G_{A,\ell}^\geo(K))$.
		We make once and for all the following remark. 
		
		\begin{rema} Let $E/K$ be a finite separable extension. Then $G_{A,\ell}(E)$ is open (and thus of finite index) in $G_{A,\ell}(K)$ and
			$G_{A,\ell}(\ol{k}E)$ is open (and thus of finite index) in $G_{A,\ell}(\ol{k}K)$. 
			The algebraic group
			$G_{A,\ell}(E)^\Zar$ is of finite index in $G_{A,\ell}(K)^\Zar=\GGG_{A,\ell}$ and
			$G_{A,\ell}(\ol{k}E)^\Zar$ is of finite index in $G_{A,\ell}(\ol{k}K)^\Zar=\GGG_{A,\ell}^\geo$. 
			In particular we have $G_{A,\ell}(E)^{\Zar, \circ}=\GGG_{A,\ell}^\circ$ and $G_{A,\ell}(\ol{k}E)^{\Zar, \circ}=\GGG_{A,\ell}^{\geo, \circ}$. It follows that 
			$\mfg_{A,\ell}=\Lie(G_{A,\ell}(E))=\Lie(G_{A,\ell}(K_A^\conn))$, $\mfg_{A,\ell}^\geo=\Lie(G_{A,\ell}^\geo(E))=\Lie(G_{A,\ell}^\geo(K_A^\conn))$,
			$\Lie(\GGG_{A,\ell})=\Lie(\GGG_{A,\ell}^\circ)$ and $\Lie(\GGG_{A,\ell}^\geo)=\Lie(\GGG_{A,\ell}^{\geo, \circ})$.
		\end{rema}
		
		\begin{lemm} \label{lemm:ZarClosure} Let $\ell\in\Ll$ and let $G$ be a compact subgroup of $\GL_{V_\ell(A)}(\Qq_\ell)$ that acts semisimply on $V_\ell(A)$. Let $\GGG=G^\Zar$ and $\mfg$ the Lie algebra of $G$.
			\begin{enumerate}
				\item[(a)] We have $\mfg^\alg=\Lie(\GGG)$. 
				\item[(b)] We have $\Lie(\GGG)=\Lie(\GGG(\Qq_\ell))$.
				\item[(c)] The group $G$ is open in $\GGG(\Qq_\ell)$ if and only if $\mf{g}=\mf{g}^\alg$. 
			\end{enumerate}
		\end{lemm} 
		
		\begin{proof} Part (a) holds true by \cite{serr-dris}*{Proposition~2~(i)} and Part (b) follows by \cite{platonov-rapinchuk}*{\S~3.1,~Lemma~3.5~and~p.~135}. If $G$ is open in $\GGG(\Qq_\ell)$, then
			$$\mfg=\Lie(G)=\Lie(\GGG(\Qq_\ell))\buildrel (b)\over = \Lie(\GGG) \buildrel (a)\over = \mfg^\alg.$$
			If $\mfg=\mfg^\alg$, then $\Lie(G)=\Lie(\GGG(\Qq_\ell))$ by Parts (a) and (b), and, as $G\subset \GGG(\Qq_\ell))$, \cite{serre-lie}*{pp.~130-131,~II.2.,~Corollary~1} then gives that $G$ is open in $\GGG(\Qq_\ell)$.
		\end{proof}
		\begin{rema}\label{remark:serre} 
			\begin{enumerate}
				\item[(a)] By Fact \ref{remark:semisimple}  the groups $G_{A, \ell}(K)$ and $G_{A, \ell}^\geo(K)$ act semisimply on $V_\ell(A)$. The Lie algebra $\mfg_{A, \ell}$ of $G_{A,\ell}(K)$ is then reductive \cite{serr-dris}*{Proposition~1}. Following \cite{serr-dris}*{pp.~119-120}, we have the direct sum decomposition $\mfg_{A,\ell}=\mfz_{A, \ell}\oplus \mfs_{A,\ell}$, where $\mfz_{A, \ell}$ is the center of $\mfg_{A,\ell}$ and $\mfs_{A,\ell}=[\mfg_{A,\ell}, \mfg_{A,\ell}]$ is semisimple and algebraic.  
				\item[(b)] By \cite{serr-dris}*{pp.~119-120} $\mfz^\alg_{A, \ell}$ is commutative and acts semisimply on $V_\ell(A)$, and thus $\mfg^\alg_{A,\ell}$ is reductive and splits up as a sum of ideals $\mfz_{A,\ell}^\alg\oplus [\mfg_{A,\ell}^\alg, \mfg_{A,\ell}^\alg]$ where $\mfz_{A,\ell}^\alg$ is the center of $\mfg_{A,\ell}^\alg$ and by Part (a), $\mfs_{A,\ell}=[\mfg_{A,\ell}^\alg, \mfg_{A,\ell}^\alg]$.
				
				\item[(c)] Lemma  \ref{lemm:ZarClosure} implies $\Lie(\GGG_{A, \ell}^\circ)=\Lie(\GGG_{A,\ell})=\mfg_{A,\ell}^\alg$ and, as $\GGG_{A,\ell}^{\geo, \circ}$ is the derived group
				of $\GGG_{A,\ell}^\circ$, it follows from 
				\cite{DG}*{II,~\S~6,~no.~2,~Propositions~2.1~and~2.3} that $\Lie(\GGG_{A,\ell}^\geo)=\Lie(\GGG_{A,\ell}^{\geo, \circ})=\mfs_{A,\ell}$ and $\Lie(\ZZZ_{A,\ell})=\mfz_{A,\ell}^\alg$. 
				\item[(d)] The algebra $\mfg_{A,\ell}$ is algebraic if and only if  $G_{A,\ell}(K)$ is $\ell$-adically open in $\GGG_{A,\ell}(\Qq_\ell)$ by Lemma \ref{lemm:ZarClosure}.
			\end{enumerate}
		\end{rema}
		
		\begin{lemm}\label{lemm:geo_open}
			\begin{enumerate}
				\item[(a)] We have $\mfg_{A,\ell}^\geo=\mfs_{A,\ell}$. 
				\item[(b)] The group $G_{A,\ell}^\geo(K)$ is open in 
				$\GGG_{A,\ell}^{\geo}(\Qq_\ell)$. 
			\end{enumerate} 
		\end{lemm}
		
		\begin{proof}
			The group $G_{A,\ell}^\geo(K)$ acts semisimply on $V_\ell(A)$ by Fact \ref{remark:semisimple} , the Zariski closure of $G_{A,\ell}^\geo(K)$ is $\GGG_{A,\ell}^\geo$ and  $\GGG_{A,\ell}^{\geo,\circ}$ is semisimple. Hence $(\mfg_{A,\ell}^\geo)^\alg=\Lie(\GGG_{A,\ell}^{\geo,\circ})$ by Lemma \ref{lemm:ZarClosure}, and that Lie algebra is semisimple. Furthermore $\mfg_{A,\ell}^\geo$ is an ideal in $(\mfg_{A,\ell}^\geo)^\alg$ \cite{DG}*{II,~\S~6,~no.~2,~Proposition~2.6~(a)}. It follows that $\mfg_{A,\ell}^\geo$ is semisimple, and thus algebraic by \cite{DG}*{II,~\S~6,~no.~2,~Corollary 2.7}. Hence
			$\mfg_{A,\ell}^\geo=(\mfg_{A,\ell}^\geo)^\alg=\Lie(\GGG_{A,\ell}^{\geo,\circ})$ and we know from Remark \ref{remark:serre} that $\Lie(\GGG_{A,\ell}^{\geo,\circ})=\mfs_{A,\ell}$. This establishes Part (a). Part (b) now follows by Lemma \ref{lemm:ZarClosure}.
		\end{proof}
		
		\begin{fact} \label{fact:dim1} We have $\dim_{\Qq_\ell}(\mfz_{A,\ell})=1$ by \cite{zar77}*{Theorem~2} and \cite{zar-finchar14}*{Section~6}. 
		\end{fact}
		
		\section{Unbounded index-subgroup function}
		We keep the notation from Section \ref{Section:weil_and_cyclotomic}. In this section, we prove that the functions 	
		$$ \Ll\to \Nn^*: \ell\mapsto \left[\HAl: \HAlq\right]$$ and $$\Ll\to \Nn^*, \ell\mapsto [\GGG_{A, \ell}(\Zz_\ell): \GGG_{A, \ell}^{(q)}]$$ are both unbounded. We thank Arno Fehm for providing the proof of the following lemma. 
		
		\begin{lemm}\label{Lemma:Not_Bounded} 
			Let $q$ be a power of $p$. 
			For $\ell\in\Ll$, let $\oq_\ell$ be its image in $\Ff_\ell^\times$ and $\langle \ol{q}_\ell\rangle$ the group it generates. The function
			$\Ll\to \Nn^* : \ell\mapsto [\Ff_\ell^\times: \langle \ol{q}_\ell\rangle]$
			is unbounded. 
		\end{lemm}
		
		\begin{proof} Let $n\in\Nn^*$, let $F_n=\Qq(\exp(2\pi i/n), q^{1/n})$. By Chebotarev's Density Theorem, there is an $\ell\in\Ll\setminus\{2\}$ such that the  prime ideal $\ell\Zz$ of $\Zz$ splits completely in $F_n$ and thus also in  $\Qq(\exp(2\pi i/n))$. Hence, $\ell \equiv 1\ \pmod{n}$ \cite{neukirch1}*{Chapter I, (10.4) Corollary} and so $[\Ff_\ell^\times:\Ff_\ell^{\times n}]=n$. The residue field of a place $v$ of $F_n$ above $\ell \Zz$ is $\Ff_\ell$ and since $q\in\Zz_\ell^\times$, its image $\ol{q}_\ell$ is an $n$-th power and so $[\Ff_\ell^\times: \langle \ol{q}_\ell\rangle]\ge [\Ff_\ell^\times:\Ff_\ell^{\times n}]$.
		\end{proof}
		
		\begin{rema} \label{rema:Zl} For odd $\ell$ we have $\Zz_\ell^\times\cong \Ff_\ell^\times \times (\Zz_\ell, +)$, and $\Zz_2^\times\cong C_2\times (\Zz_2,+)$. Hence, for every $\ell\in\Ll$, every infinite closed subgroup of $\Zz_\ell^\times$ is open in $\Zz_\ell^\times$ and every closed subgroup of $\Zz_\ell^\times$ can be topologically generated by two elements. 
		\end{rema}
		\begin{prop}\label{Prop:Unbounded_Homothety_Function} For every $\ell \in \Ll$, $\langle q_\ell\rangle$ is open in $\Zz_\ell^\times$ and $\HAlq $ is open in $\HAl$. The function $ \Ll\to \Nn^*: \ell\mapsto \left[\HAl: \HAlq\right]$ is unbounded. 
		\end{prop}
		
		\begin{proof} The prime power $q$ is not a root of unity in $\Zz_\ell^\times$, thus $\langle q_\ell\rangle$ is infinite, hence open in $\Zz_\ell^\times$ by Remark \ref{rema:Zl}. By Lemma \ref{lemm:mdet}, the map $\diag_{A,\ell}: \Zz_\ell^\times\to \HAl,\ a\mapsto a\Id_{V_\ell(A)}$ is a group isomorphism that sends  $S_\ell=\{a\in\Zz_\ell^\times: a^2\in \langle q_\ell \rangle\}$ to $\HAlq$. Thus, $\left[\HAl: \HAlq\right]=[\Zz_\ell^\times:S_\ell]$. The group $S_\ell$ contains $q_\ell$, hence is open in $\Zz_\ell^\times$. It follows that $\HAlq$ is open in $\HAl$. The Abelian group $S_\ell/\langle q_\ell \rangle$ is $2$-torsion and is topologically generated by two elements. Hence $[S_\ell:\langle q_\ell \rangle]\le 4$. It follows that 
			$$
			\left[\HAl: \HAlq\right]=[\Zz_\ell^\times:S_\ell]\ge \frac{1}{4}[\Zz_\ell^\times:\langle q_\ell \rangle]
			\ge \frac{1}{4} [\Ff_\ell^\times: \langle \oq_\ell\rangle],
			$$
			so the assertion follows from Lemma \ref{Lemma:Not_Bounded}.
		\end{proof}
		
		\begin{lemm} \label{lemm:CenterHomotheties} Let $\ell\in\Ll$. 
			\begin{enumerate} 
				\item[(a)] The algebraic group  $\ZZZ_{A, \ell}^\circ$ contains $\HHH_{A,\ell}$ as closed subgroup scheme.
				\item[(b)] The center $\mfz_{A,\ell}^\alg$ of $\mfg_{A,\ell}^\alg$ contains $\Qq_\ell \Id_{V_\ell(A)}$. 
			\end{enumerate}
		\end{lemm}
		
		\begin{proof} We use the notation from Fact \ref{fact:compatible}. Let $s$ be a closed point of $S$ and $\TTT_\ell$ the Zariski closure of $\langle \Fr_{A,s,\ell}\rangle$ in $\GGG_{A,\ell}$. 
			The representation of $\Gal(k(s))$ on $V_\ell(A)$ is isomorphic to the representation of $\Gal(k(s))$ on $V_\ell(\AAA_s)$. It follows by Fact \ref{remark:semisimple} applied to 
			$\AAA_s/k(s)$ that the action of $\langle \Fr_{A,s,\ell}\rangle$ on $V_\ell(A)$ is semisimple and $\TTT_\ell^\circ$ is reductive. Moreover $\TTT_\ell$ is commutative. Thus $\TTT_\ell^\circ$ is a torus - the Frobenius torus at $s$. Let $m=[\TTT_\ell(\Qq_\ell):\TTT_\ell^\circ(\Qq_\ell)]$. Then $\TTT_\ell^\circ$ is the Zariski closure of $\langle \Fr_{A,s,\ell}^m\rangle$. We choose a basis of $V_\ell(A)$ and identify $\GL_{V_\ell(A)}$ with $\GL_{2g_A, \Qq_\ell}$ in this proof.  
			There exists $a_1,\cdots, a_{2g_A}\in\ol{\Qq}_\ell$ such that in  $\GL_{2g_A}(\ol{\Qq}_\ell)$ the matrix $\Fr_{A,s,\ell}$ is conjugate to $D:=\mathrm{diag}(a_1,\cdots, a_{2g_A})$. The $a_j$ are integral over $\Zz$ and for every embedding $i: \Qq(a_1,\cdots, a_{2g_A})\to \Cc$ we have $|i(a_j)|=(\# k(s))^{\frac{1}{2}}$. Thus $\TTT_{\ell, \ol{\Qq}_\ell}^\circ$ is conjugate in $\GL_{2g_A, \ol{\Qq}_\ell}$ to $\DDD:=\langle D^m\rangle^\Zar$.
			
			We claim that $\DDD$ contains $\HHH_{A,\ell, \ol{\Qq}_\ell}$. Let $\SSS$ be a copy of $\Gg_{m, \ol{\Qq}_\ell}^{2g_A}$ embedded diagonally in $\GL_{2g_A, \ol{\Qq_\ell}}$. It is enough to prove that every character $\chi: \SSS\to \Gg_{m, \ol{\Qq}_\ell}$ that vanishes on $\DDD$ vanishes on $\HHH_{A,\ell, \ol{\Qq}_\ell}$.  Let $\chi$ be such a character. 
			Let $\chi_j: \SSS\to \Gg_{m, \ol{\Qq}_\ell}$ be the projection on the $j$-th factor. There exist $n_j\in\Zz$ such that $\chi=\prod_{j=1}^{2g_A} \chi_j^{n_j}$. The character $\chi$ vanishes on $\DDD$. Hence $1=\chi(D^m)=\prod_{j=1}^{2g_A} a_j^{mn_j}$. It follows that 
			$$
			1=\prod_{j=1}^{2g_A} |i(a_j)|^{mn_j}=\prod_{j=1}^{2g_A} (\#k(s))^{\frac{1}{2}mn_j}=(\#k(s))^{\frac{1}{2} m \sum_{j=1}^{2g_A} n_j}.
			$$
			This forces $\sum_{j=1}^{2g_A} n_j=0$. Thus $\chi$ annihilates $\HHH_{A,\ell, \ol{\Qq}_\ell}$, establishing the claim. 
			
			We now know that $\DDD$ contains the normal subgroup $\HHH_{A,\ell, \ol{\Qq}_\ell}$ of $\GL_{2g_A, \ol{\Qq}_\ell}$ and moreover $\TTT_{\ell, \ol{\Qq}_\ell}^\circ$ is conjugate in $\GL_{2g_A, \ol{\Qq}_\ell}$ 
			to 
			$\DDD$. This implies that $\HHH_{A,\ell, \ol{\Qq}_\ell}$ is a closed subscheme $\TTT_{\ell, \ol{\Qq}_\ell}^\circ$. Therefore, by Lemma \ref{Technical_Lemma_Inclusion} below $\HHH_{A,\ell}$ is a closed subscheme of $\TTT_{\ell}^\circ$, and $\TTT_{\ell}^\circ$ is a closed subscheme of $\GGG_{A,\ell}^\circ$.  As $\HHH_{A,\ell}$ is central even in $\GL_{{2g_A}, \Qq_\ell}$ it follows that $\HHH_{A,\ell}\subset \ZZZ_{A,\ell}$. As $\HHH_{A,\ell}$ is connected we in fact get $\HHH_{A,\ell}\subset \ZZZ_{A,\ell}^\circ$. This proves Part (a). 
			
			From $\HHH_{A,\ell}\subset \ZZZ_{A,\ell}$ it follows by \cite{DG}*{II,~\S~6,~no.~2,~Propositions~2.1~and~2.3} that 
			$\Qq_\ell \Id_{V_\ell(A)}=\Lie(\HHH_{A,\ell})\subset \Lie( \ZZZ_{A,\ell})$ and $\Lie( \ZZZ_{A,\ell})=\mfz_{A,\ell}^\alg$ by Part (c) of Remark \ref{remark:serre}. Hence Part (b) holds. 
		\end{proof}
		
		\begin{lemm}\label{Technical_Lemma_Inclusion}
			Let $L$ be a field  and let $X$ be a $L$-variety. Let $Y$ and $Z$ be closed reduced $L$-subschemes of $X$. Let $F/L$ be a field extension and $\Omega/F$ an algebraically closed extension. The following are equivalent. 
			\begin{enumerate}
				\item[(i)] $Y$ is a closed subscheme of $Z$. 
				\item[(ii)] $Y(\Omega)\subset Z(\Omega)$ 
				\item[(iii)] $Y_F$ is a closed subscheme of $Z_F$. 
			\end{enumerate}
		\end{lemm}

		\begin{proof} We identify $Y_F$ and $Z_F$ with closed subschemes of $X_F$. Implication (i)$\Rightarrow$(ii) is by functor of points. Now assume (ii) and let $y$ be a closed point of $Y$. Choose an embedding of the residue field $\kappa(y)$ into $\Omega$. We get a geometric point $y': \Spec(\Omega)\to \Spec(\kappa(y))\to Y\to X$ localized at $y$, and by (ii) $y'$ factors through $Z$. Hence $y\in Z$. We have shown that every closed point of  $Y$ lies in the underlying topological space of $Z$. As the set of all closed points of $Y$ is dense in $Y$ it follows that the space underlying $Y$ is contained in the space underlying $Z$. As $Y$ and $Z$ are reduced, (i) follows by \cite{grothendieck_1960}*{I.5.2.4}. The same proof strategy gives (ii)$\Leftrightarrow$(iii).  
		\end{proof}
		
		\begin{coro} \label{coro:notbounded2} The function
			$\Ll\to \Nn^*, \ell\mapsto [\GGG_{A, \ell}(\Zz_\ell): \GGG_{A, \ell}^{(q)}] $ is unbounded.
		\end{coro} 
		
		\begin{proof} We have $\HHH_{A, \ell}(\Zz_\ell)\subset \GGG_{A,\ell}(\Zz_\ell)$ by Lemma \ref{lemm:CenterHomotheties} and 
			$\HHH_{A,\ell}^{(q)}=\HHH_{A, \ell}(\Zz_\ell)\cap \GGG_{A,\ell}^{(q)}$. Thus 
			$ [\GGG_{A, \ell}(\Zz_\ell): \GGG_{A, \ell}^{(q)}] \ge [\HHH_{A, \ell}(\Zz_\ell): \HHH_{A, \ell}^{(q)}]$ and readily the latter index is unbounded by Proposition 
			\ref{Prop:Unbounded_Homothety_Function}. 
		\end{proof}
		
		\begin{coro} \label{coro:notbounded} 
			\begin{enumerate}
				\item[(a)] The function $\Ll\to \Nn^* \cup \{\infty\}:\ \ell \mapsto [\HAl:\CAl(K)]$  is unbounded.
				\item[(b)] The function $\Ll\to \Nn^* \cup \{\infty\}:\ \ell \mapsto [\GGG_{A,\ell}(\Zz_\ell): G_{A,\ell}(K)]$  is unbounded.
			\end{enumerate}
		\end{coro}
		
		\begin{proof} By Corollary \ref{rema:topo_gen}, we have $\CAl(K) \subset \HAlq\subset  \HAl$ and
			$G_{A,\ell}(K)\subset \GGG_{A,\ell}^{(q)}\subset \GGG_{A,\ell}(\Zz_\ell)$. Thus Part (a) follows by Proposition \ref{Prop:Unbounded_Homothety_Function} and
			Part (b) by Corollary \ref{coro:notbounded2}. 
		\end{proof}
		
		\section{Homothety condition and independence} We keep the notation from Section \ref{Section:weil_and_cyclotomic}. In this section we offer various equivalent characterizations of condition \condH.
		
		\begin{lemm}  \label{lemm:Hl} Let $\ell\in\Ll$. The following conditions are equivalent. 
			\begin{enumerate}
				\item[(a)] The Lie algebra $\mfg_{A,\ell}$ contains $\Qq_\ell \cdot \Id_{V_\ell(A)}$. 
				\item[(b)] The center $\mfz_{A,\ell}$ of $\mfg_{A,\ell}$ equals $\Qq_\ell \cdot \Id_{V_\ell(A)}$.
				\item[(c)] We have $\dim_{\Qq_\ell}(\mfz_{A,\ell}^\alg)=1$. 
				\item[(d)] The Lie algebra $\mfz_{A,\ell}$ is algebraic. 
				\item[(e)] The Lie algebra $\mfg_{A,\ell}$ is algebraic. 
				\item[(f)] We have $\dim(\ZZZ_{A,\ell})=1$. 
				\item[(g)] The identity component $\ZZZ_{A,\ell}^\circ$ of $\ZZZ_{A,\ell}$ equals $\HHH_{A,\ell}$.
				\item[(h)] The group $\CAl(K)$ is open in $\HAl$. 
				\item[(i)] The group $\CAl(K)$ is infinite.
				\item[(j)] The group $G_{A,\ell}(K)$ is open in $\GGG_{A,\ell}(\Zz_\ell)$. 
			\end{enumerate}
		\end{lemm}
		
		\begin{proof}
			$(a)\Rightarrow(b)$: 
			Since (a) holds, then $\mfz_{A,\ell}$ contains $\Qq_\ell \cdot \Id_{V_\ell(A)}$. But, 
			by Fact \ref{fact:dim1}, $\mfz_{A,\ell}$ is one-dimensional over $\Qq_\ell$. Hence (b) holds. 
			
			$(b)\Rightarrow(a)$: 
			Since (b) holds, then $\Qq_\ell\Id_{V_\ell(A)}=\mfz_{A,\ell}\subset \mfg_{A,\ell}$. Hence (a) holds.
			
			$(b)\Leftrightarrow (c)\Leftrightarrow (d)$ We know that $\mfz_{A,\ell}\subset \mfz_{A,\ell}^\alg$, and $\mfz_{A,\ell}^\alg\supset \Qq_\ell \Id_{V_\ell(A)}$ by Lemma
			\ref{lemm:CenterHomotheties}. Moreover $\dim(\mfz_{A,\ell})=1$ by Fact \ref{fact:dim1}, and $\dim_{\Qq_\ell}( \Qq_\ell \Id_{V_\ell(A)})=1$. Hence 
			$$
			\mfz_{A,\ell}^\alg=\Qq_\ell \Id_{V_\ell(A)}\iff \dim(\mfz_{A,\ell}^\alg)=1\ \iff \mfz_{A,\ell}=\mfz_{A,\ell}^\alg.
			$$
			
			$(d)\Leftrightarrow (e)$. This equivalence is immediate from Remark \ref{remark:serre}. 
			
			$(c)\Leftrightarrow(f)$: By Remark \ref{remark:serre} we have $\Lie(\ZZZ_{A,\ell})=\mfz_{A,\ell}^\alg$. This implies
			$\dim(\ZZZ_{A,\ell})=\dim_{\Qq_\ell}(\mfz_{A,\ell}^\alg)$. Thus we have $\dim_{\Qq_\ell}(\mfz_{A,\ell}^\alg)=1\ \iff \dim(\ZZZ_{A,\ell})=1$, as desired for the 
			equivalence  $(c)\Leftrightarrow(f)$. 
			
			$(f)\Leftrightarrow(g)$: By Lemma \ref{lemm:CenterHomotheties} there is the inclusion $\HHH_{A,\ell}\subset \ZZZ_{A,\ell}^\circ$, and $\dim(\HHH_{A,\ell})=1$. 
			Moreover $ \dim(\ZZZ_{A,\ell}^\circ)=\dim(\ZZZ_{A,\ell})$. 
			Thus
			$\HHH_{A,\ell}= \ZZZ_{A,\ell}^\circ\ \iff \ \dim(\ZZZ_{A,\ell}^\circ)=1\ \iff \ \dim(\ZZZ_{A,\ell})=1$, as desired. 
			
			$(a)\Leftrightarrow (h)$:  The Lie algebra $\mfc_{A,\ell}$ of $\CAl(K)$ verifies $\mfc_{A,\ell}=\mfg_{A,\ell}\cap \Qq_\ell\cdot\Id_{V_\ell(A)}$ by
			\cite{serre-lie}*{p.~130,~II.2,~Theorem}. Thus $(a)\ \iff\ \mfc_{A,\ell} = \Qq_\ell \cdot \Id_{V_\ell(A)}$. 
			Moreover $\CAl(K)\subset \HAl$. Hence 
			\cite{serre-lie}*{pp.~130-131,~II.2.,~Corollary~1} gives $\mfc_{A,\ell} = \Qq_\ell \cdot \Id_{V_\ell(A)}\ \iff (h)$. 
			
			$(h)\Leftrightarrow (i)$: The profinite group $\HHH_{A,\ell}(\Zz_\ell)$ is isomorphic to $\Zz_\ell^\times$. It follows from Remark \ref{rema:Zl} that the closed subgroup $C_{A, \ell}(K)$ of $\HHH_{A,\ell}(\Zz_\ell)$ is open in $\HHH_{A,\ell}(\Zz_\ell)$ if and only if it is infinite. 
			
			$(e)\Leftrightarrow (j)$: This equivalence holds by Lemma \ref{lemm:ZarClosure}.		
		\end{proof}
		
		Our next goal, Lemma \ref{lemm:H2}, is to prove that if the equivalent conditions of Lemma \ref{lemm:Hl} hold for one $\ell \in \Ll$, they hold for all $\ell \in \Ll$. We start with intermediate results.

		The following theorem is due to Chin who proved it for curves in \cite{chin_2004}*{Theorem 1.4 (ii)} and reduced it to that case by a Katz's space-filling curves argument \cite{chin_2004}*{p.~726}. For the reader's convenience, we include an alternative reduction to curves. 
		
		\begin{thm} {\bf (Chin)} \label{thm:chin} There is a number field $F$, an $F$-vector space $V(A)$, and a split reductive subgroup $\GGG_A$ of $\GL_{V(A)}$ and, for 
			every $\ell\in\Ll$ and every  place $\lambda$ of $F$ over $\ell$, an isomorphism $(I_\lambda, f_\lambda)$ of faithful representations
			$$
			(\GGG_A\to \GL_{V(A)})\times_F F_\lambda\ \cong\ (\GGG_{A,\ell}^\circ\to \GL_{V_\ell(A)})\times_{\Qq_\ell} F_\lambda.
			$$
			In other words, we have algebraic group isomorphisms $I_\lambda: \GGG_A\times_F F_\lambda\to \GGG_{A, \ell}^\circ \times_{\Qq_\ell} F_\lambda$ and $F_\lambda$-vector space isomorphisms
			$f_\lambda: V(A)\otimes_F F_\lambda\to V_\ell(A)\otimes_{\Qq_\ell} F_\lambda$ such that the diagram
			
			$$
			\xymatrix{
				\GGG_A\times_F F_\lambda \ar[r]^{}\ar[d]^{I_\lambda} & \GL_{V(A)\otimes_{F} F_\lambda}\ar[d]^{\GL_{f_\lambda}}\\
				\GGG_{A, \ell}^\circ\times_{\Qq_\ell} F_\lambda \ar[r]^{} & \GL_{V_\ell(A)\otimes_{\Qq_\ell} F_\lambda}.
			}$$
			commutes. 
		\end{thm}
		
		\begin{proof} By \cite{bgp2}*{Corollary 3.33} there is a connected finite \'etale cover $g: S'\to S$, a smooth connected curve $C$ and an embedding $h: C\to S'$ such that for every $\ell\in\Ll$ 
			the image of
			$$
			\pi_1(C)\buildrel h_*\over\longrightarrow 
			\pi_1(S') \buildrel g_*\over\longrightarrow   \pi_1(S)
			\buildrel \rho_{A,\ell} \over\longrightarrow  \GL_{V_\ell(A)}(\Qq_\ell)
			$$
			equals $\rho_{A,\ell}(g_*(\pi_1(S')))$. Now, the group $g_*(\pi_1(S'))$ is open in $\pi_1(S)$. It follows that $\rho_{A,\ell}\circ g_*\circ h_*(\pi_1(C))$ is open in $G_{A,\ell}(K) = \rho_{A,\ell}(\pi_1(S))$. Hence
			$((\rho_{A,\ell}\circ g_*\circ h_*(\pi_1(C)))^{\Zar})^\circ=\GGG_{A,\ell}^\circ$. 
			Let $(\LL_\ell)_{\ell\in \Ll}$ be the $\Qq$-compatible system of Remark \ref{rema:chin_applies}. Then $(h^*g^*\LL_\ell)_{\ell\in\Ll}$ is a pure semisimple $\Qq$-compatible system and $\rho_{A, \ell}\circ g_*\circ h_*$ is the representation of $\pi_1(C)$ associated to the lisse $\Qq_\ell$-sheaf $h^*g^* \LL_\ell$. The assertion follows by \cite{chin_2004}*{Theorem 1.4 (ii)}.
		\end{proof}
		
		\begin{defi} Let $\ZZZ_A$ be the center of the reductive group $\GGG_A$ from Theorem \ref{thm:chin}. 
		\end{defi}
		
		\begin{rema} \label{rema:centerchin} Let $\ell\in \Ll$. 
			For an algebraic group over a field of characteristic zero the formation of the center commutes with base change. Hence $\dim(\ZZZ_A)=\dim(\ZZZ_{A,\ell})$.  It follows
			that the equivalent conditions (a)-(j) from Lemma \ref{lemm:Hl} are also equivalent to the following condition 
			\begin{enumerate}
				\item[(k)] $\dim(\ZZZ_A)=1$. 
			\end{enumerate}
		\end{rema}
		
		\begin{lemm} \label{lemm:H2} The following are equivalent.
			\begin{enumerate}
				\item[(a)] The equivalent conditions of Lemma \ref{lemm:Hl} hold for some $\ell \in \Ll$.
				\item[(b)] The equivalent conditions of Lemma \ref{lemm:Hl} hold for all $\ell \in \Ll$.
				\item[(c)] The Abelian variety $A$ satisfies condition \condH. 
			\end{enumerate}
		\end{lemm}
		
		\begin{proof}
			The equivalence $(a)\Leftrightarrow (b)$ holds because the equivalent conditions (a)-(j) of Lemma \ref{lemm:Hl} are equivalent to
			condition (k) of Remark \ref{rema:centerchin}, and the latter does not depend on $\ell\in\Ll$. 
			Moreover condition (b) of Lemma \ref{lemm:Hl}  holds for all $\ell\in\Ll$ if and only if $A$ satisfies condition \condH. Hence $(b)\Leftrightarrow (c)$
		\end{proof}
		
		\section{Proof of the main results}
		
		\begin{defi}	
			For $\ell \in \Ll$, let $\KKK_{A,\ell}$ be the kernel of the restriction $m_{A,\ell}|\GGG_{A,\ell}^\circ: \GGG_{A, \ell}^\circ\to \Gg_{m,\Qq_\ell}$.
		\end{defi}
		
		\begin{lemm} \label{lemm:boundK}
			\begin{enumerate}
				\item[(a)] The group scheme $\KKK_{A,\ell}$ contains $\GGG_{A,\ell}^{\geo, \circ}$ as a normal subgroup scheme.
				\item[(b)] If $A/K$ satisfies condition \condH, then the group scheme  $\KKK_{A,\ell}/\GGG_{A,\ell}^{\geo, \circ}$ is finite \'etale  over $\Qq_\ell$ and 
				$\rk_{\Qq_\ell}(\KKK_{A,\ell}/\GGG_{A,\ell}^{\geo, \circ})\le 2$.
			\end{enumerate}
		\end{lemm}
		
		\begin{proof}
		Let $\chi_\ell$ be the restriction $m_{A,\ell}|\GGG_{A,\ell}^\circ: \GGG_{A, \ell}^\circ\to \Gg_m$ of $m_{A, \ell}$ to the reductive group 
			$\GGG_{A, \ell}^\circ$. The algebraic group $\GGG_{A,\ell}^{\geo, \circ}$ is the derived subgroup of the reductive group $\GGG_{A,\ell}^\circ$ (Fact \ref{remark:derived}). 
			In particular $\GGG_{A,\ell}^{\geo, \circ}$ is semisimple. Hence $\chi_\ell$ annihilates $\GGG_{A,\ell}^{\geo, \circ}$. Part (a) thus follows. 
			
			Assume from now on that $A/K$ satisfies condition \condH. Let $\GGG^\ab=\GGG_{A,\ell}^\circ/ \GGG_{A,\ell}^{\geo, \circ}$. The character $\chi_\ell$  factors as 
			$$\GGG_{A,\ell}^{\circ} \buildrel \pi\over \longrightarrow \GGG^\ab\buildrel \ol{\chi}_\ell \over
			\longrightarrow \Gg_m,$$
			
			and $\KKK_{A,\ell}/\GGG_{A,\ell}^{\geo, \circ}\cong \ker(\ol{\chi}_\ell)$. It is thus enough to prove that  $\ker(\ol{\chi}_\ell)$ is finite \'etale over $\Qq_\ell$ of rank $\le 2$. 
			Let $\pi'$ be the restriction $\pi|\ZZZ_{A,\ell}^\circ: \ZZZ_{A,\ell}^\circ\to  \GGG^\ab$. Now, $\GGG_{A,\ell}^\circ = \ZZZ_{A,\ell}^\circ \cdot \GGG_{A,\ell}^{\geo, \circ}$ by \cite{Borel_1991}*{p.~181,~14.2~Proposition~(3)}. Hence $\pi'$ is surjective. It follows $\pi'$ restricts to a surjective homomorphism of group schemes 
			$\ker(\ol{\chi_\ell}\circ \pi')\to \ker(\ol{\chi}_\ell)$. By condition \condH~we have $\ZZZ_{A,\ell}^\circ=\HHH_{A,\ell}$, and 
			$\ol{\chi_\ell}\circ \pi'(x\Id_{V_\ell(A)})=\chi_\ell(x \Id_{V_\ell(A)})=x^2$, hence $\ker(\ol{\chi_\ell}\circ \pi')\cong \mu_{2, \Qq_\ell}$. It follows that 
			$\ker(\ol{\chi}_\ell)\cong\mu_{2, \Qq_\ell}$ or $\ker(\ol{\chi}_\ell)$ is trivial, as desired. 
		\end{proof}
		
		\begin{defi} If $G$ is a closed subgroup of $\GL_{T_\ell(A)}(\Zz_\ell)$ we denote by $G^+$ the closed subgroup of $G$ topologically generated by the $\ell$-Sylow subgroups of the profinite group $G$. 
		\end{defi}
		
		\begin{fact}\label{fact:zarhin-mod-l} For all but finitely many $\ell\in\Ll$ the $\Ff_\ell[\Gal(K)]$-module $A[\ell](K^\sep)$ is semisimple by  \cite{zar-finchar14}*{Corollary 2.3 (iii)}. 
		\end{fact}
		
		\begin{fact}\label{fact:CHT} 
			There is an integer $C_{2g_A} \geq 1$ depending only on $g_A$, a finite separable extension $K'_A/K$ and a finite subset $\Sigma_0$ of $\Ll$ such that the following statements hold. 
			\begin{enumerate}
				\item[(a)] The Zariski closure $\mfG_{A,\ell}^\geo$ of $G_{A,\ell}^\geo(K'_A)$ in $\GL_{T_\ell(A)}$ is a semisimple $\Zz_\ell$-group scheme for all $\ell\in\Ll\setminus \Sigma_0$ by \cite{cht}*{Corollary 7.5} and Fact \ref{fact:zarhin-mod-l}. 
				\item[(b)] We have $[\mfG_{A,\ell}^\geo(\Zz_\ell): G_{A,\ell}^\geo(K'_A)]\le C_{2g_A}$ for all $\ell\in \Ll\setminus \Sigma_0$  by \cite{cht}*{Theorem 7.3 (7.3.2)}. 
				\item[(c)] We have $G_{A, \ell}^\geo(K'_A)= \mfG_{A,\ell}^\geo(\Zz_\ell)^+$ for all $\ell\in\Ll\setminus \Sigma_0$ by \cite{cht}*{Theorem 7.3 (7.3.2)}. 
			\end{enumerate}
		\end{fact}
		
		\begin{lemm} \label{lemma:CHTplusind} There exists a finite separable extension $K_A^\ind$ of $K$ and a finite subset $\Sigma=\Sigma(A/K)$ of $\Ll$ depending on $A/K$ with the following properties. 
			\begin{enumerate}
				\item[(a)] We have $G_{A,\Ll}^\geo(K_A^\ind)=\prod_{\ell\in\Ll} G_{A, \ell}^{\geo}(K_A^\ind)$. 
				\item[(b)] The group  $G_{A, \ell}^\geo(K_A^\ind)$ is a characteristic open subgroup of $\GGG_{A,\ell}^{\geo, \circ}(\Zz_\ell)$ for every $\ell\in\Ll$. 
				\item[(c)] We have $[\GGG_{A,\ell}^{\geo, \circ}(\Zz_\ell): G_{A,\ell}^\geo(K^\ind_A)]\le C_{2g_A}$ for all $\ell\in \Ll\setminus \Sigma$, where $C_{2g_A}$ is the constant from
				Fact \ref{fact:CHT}. 
				\item[(d)] We have $K_A^\ind\supset K_A^{(3,5,7)}\supset K_A^\conn$ and $\ol{k}K_A^\ind\supset K_A^{\geo, \conn}$. 
			\end{enumerate}
		\end{lemm}
		
		\begin{proof}
			By \cite{bgp}*{Theorem~7.7} (applied over $K_A' \cdot K_A^{(3,5,7)}$) there is a finite separable extension $K''/K$ containing $K_A'$ and $K_A^{(3,5,7)}$ 
			(*) such that   
			$G_{A, \Ll\smallsetminus \{2,3\}}^\geo(K'')=\prod_{\ell\in \Ll\setminus \{2,3\}} G_{A, \ell}^\geo(K'')$. 
			Up to replacing $K''$ by a finite separable extension, \cite{serre-ind}*{Lemme 3} implies that 
			$G_{A, \Ll}^\geo(K'')=\prod_{\ell\in \Ll} G_{A, \ell}^\geo(K'')$. Let $\Sigma_0$ be the subset of $\Ll$ from Fact \ref{fact:CHT}. 
			Let
			$$
			\Sigma :=\Sigma_0\cup\{\ell\in\Ll: \mbox{$\ell\le [\ok K'':\ok K_A']$}\}.
			$$
			Let $\ell\in\Ll\setminus \Sigma$. 
			Then $G_{A,\ell}^\geo(K_A')=\mfG_{A,\ell}^{\geo}(\Zz_\ell)^+$ and thus $G_{A,\ell}^\geo(K_A')=G_{A,\ell}^\geo(K_A')^+$. The index 
			$[G_{A,\ell}^\geo(K_A'):G_{A,\ell}^\geo(K'')]$ is $<\ell$, hence each $\ell$-Sylow subgroup of $G_{A,\ell}^\geo(K_A')$ acts trivially on the coset space $G_{A,\ell}^\geo(K_A')/G_{A,\ell}^\geo(K'')$. It
			follows that each $\ell$-Sylow subgroup of $G_{A,\ell}^\geo(K_A')$  is contained in 
			$G_{A,\ell}^\geo(K'')$. Hence $G_{A,\ell}^\geo(K'')=G_{A,\ell}^\geo(K_A')=\mfG_{A,\ell}^\geo(\Zz_\ell)^+$ for all 
			$\ell\in\Ll\setminus \Sigma$. The generic fibre of $\mfG_{A,\ell}^\geo/\Zz_\ell$ is the Zariski closure of $G_{A,\ell}^\geo(K_A')$ in $\GL_{V_\ell(A)}$ and $K_A'/K$ is finite. 
			Moreover,  $\mfG_{A,\ell}^\geo/\Zz_\ell$ is semisimple and thus the generic fibre of $\mfG_{A,\ell}^\geo/\Zz_\ell$ is connected. 
			Hence the generic fibre of  $\mfG_{A,\ell}^\geo/\Zz_\ell$ equals $\GGG_{A,\ell}^{\geo, \circ}$, and 
			thus $\mfG_{A,\ell}^\geo(\Zz_\ell)=\GGG_{A,\ell}^{\geo, \circ}(\Zz_\ell)$. Altogether we see that 
			$G_{A,\ell}^\geo(K'')=G_{A,\ell}^\geo(K_A')=\GGG_{A,\ell}^{\geo, \circ}(\Zz_\ell)^+$ for all 
			$\ell\in\Ll\setminus \Sigma$. 
			In particular, $G_{A,\ell}^\geo(K'')$ is a characteristic open subgroup of $\GGG_{A,\ell}^{\geo, \circ}(\Zz_\ell)$,
			and $[\GGG_{A,\ell}^{\geo, \circ}(\Zz_\ell): G_{A,\ell}^\geo(K'')]\le C_{2g_A}$ for all 
			$\ell\in\Ll\setminus \Sigma$.
			
			Let $\ell\in\Ll$. Lemma   \ref{lemm:Kconnbound} and (*) imply $\ol{k} K''\supset K_A^{\geo,\conn}$, and $\ol{k} K''|K_A^{\geo,\conn}$ is a finite extension. 
			Hence $G_{A, \ell}^\geo(K'')\subset \GGG_{A, \ell}^{\geo, \circ}(\Zz_\ell)$ and $\GGG_{A,\ell}^{\geo, \circ}=G_{A, \ell}^\geo(K'')^\Zar$. 
			Now Lemma \ref{lemm:geo_open} implies that
			$G_{A,\ell}^\geo(K'')$ is open in $\GGG_{A, \ell}^{\geo, \circ}(\Zz_\ell)$ . 
			
			For $\ell\in \Sigma$ we define
			$$N_{A, \ell}:=\bigcap_{\alpha\in\Aut_{\mathrm{cont}}(\GGG_{A, \ell}^{\geo, \circ}(\Zz_\ell))} \alpha(G_{A,\ell}^\geo(K'')).$$
			The profinite group
			$\GGG_{A, \ell}^{\geo, \circ}(\Zz_\ell)$ is topologically finitely generated (see Fact \ref{fact:fingen}) and thus has, for every positive integer $s$, only finitely many open subgroups of index $s$. Hence $N_{A,\ell}$ is open in $\GGG_{A, \ell}^{\geo, \circ}(\Zz_\ell)$ for all $\ell\in\Sigma$. By construction  $N_{A,\ell}$ is a characteristic subgroup of $\GGG_{A, \ell}^{\geo, \circ}(\Zz_\ell)$ for all $\ell\in\Sigma$. 
			
			For $\ell\in\Ll\setminus \Sigma$ we define $N_{A,\ell}=G_{A,\ell}^\geo(K'')$. Then $N_{A, \ell}$ is a characteristic open subgroup of $ \GGG_{A, \ell}^{\geo, \circ}(\Zz_\ell)$ for every $\ell\in\Ll$ and 
			$[\GGG_{A,\ell}^{\geo, \circ}(\Zz_\ell): N_{A,\ell}]\le C_{2g_A}$ for all $\ell\in\Ll\setminus \Sigma$. The group $N_{A,\Ll}:=\prod_{\ell\in\Ll} N_{A,\ell}$ is an open subgroup of $G_{A, \Ll}^\geo(K'')=\prod_{\ell\in \Ll} G_{A, \ell}^\geo(K'')$ which is in turn an open subgroup of $G_{A,\Ll}^\geo(K)$. We let $E$ be the fixed field of $\rho_{A,\Ll}^{-1}(N_{A,\Ll})\cap\Gal(\ok K'')$ in $K^\sep$. Then $E/\ok K''$ is a finite extension and $\rho_{A,\Ll}(\Gal(E))=N_{A, \Ll}$. Choose a $\ok K''$-basis $(a_1,\cdots, a_t)$ of $E$ and let $K_A^\ind=K''(a_1,\cdots, a_t)$. Then$E=\ok K_A^\ind$ and thus $G_{A, \Ll}^\geo(K_A^\ind)=\rho_{A, \Ll}(\Gal(\ok K_A^\ind))=N_{A,\Ll}$, and the assertions (a), (b) and~(c) hold. 
			Part (d) then follows by (*) and Lemma \ref{lemm:Kconnbound}.
		\end{proof} 
		
		Fix once and for all a field $K_A^\ind$ as in Lemma \ref{lemma:CHTplusind}. 
		
		\begin{defirema} \label{defrem:constant} For every subset $L$ of $\Ll$, let 
			$$
			c(A/K, L):=\sup\{[\GGG_{A, \ell}^{\geo, \circ}(\Zz_\ell):G_{A,\ell}^\geo(K_A^\ind)]: \ell\in L \}.
			$$
			Parts (b) and (c) of Lemma \ref{lemma:CHTplusind}  imply that $c(A/K, L)<\infty$. If $L\subset \Ll\setminus \Sigma$, then $c(A/K, L)\le C_{2g_A}$. 
		\end{defirema}
		
		\begin{lemm}\label{lemm:grp:index}
			Let $M\subset N$ be normal subgroups of a group $G$ and assume that
			$[N:M]<\infty.$ Let $v\in\Nn^*$ with $v\ge [N:M]$. 
			Then, for every $g\in G$ and every $n\in N$, one has
			$(gn)^{v!}\in g^{v!}M$. 
		\end{lemm}
		
		\begin{proof} Let $u=[N:M]$. 
			Since $M$ and $N$ are normal in $G$, conjugation by $g$ induces an
			automorphism
			$\alpha:N/M\rightarrow N/M,\ 
			xM\mapsto gxg^{-1}M.$
			Consider the map
			$
			T:N/M\rightarrow N/M,\ 
			xM\mapsto \alpha(nM)\alpha(xM).
			$
			This is a permutation of $N/M$, since $\alpha$ is an automorphism and
			left multiplication by $\alpha(nM)$ is bijective.
			
			The set $N/M$ has cardinality $u$. Hence
			$T^{u!}=\operatorname{id}_{N/M}.$
			Put $r=u!$. Iterating $T$ and evaluating at the neutral coset $M$ gives
			$$M=T^r(M)
			=\alpha(nM)\alpha^2(nM)\cdots\alpha^r(nM).$$
			Since
			$\alpha^j(nM)=g^jng^{-j}M,$
			we obtain
			$
			(gng^{-1})(g^2ng^{-2})\cdots(g^rng^{-r})\in M.
			$
			On the other hand,
			$(gn)^rg^{-r}
			=
			(gng^{-1})(g^2ng^{-2})\cdots(g^rng^{-r}).$
			Therefore
			$(gn)^rg^{-r}\in M.$
			As $M$ is normal in $G$, this is equivalent to
			$(gn)^r\in g^rM.$ Thus there exists $m\in M$ with $(gn)^r=g^rm$. 
			From $v\ge u$ and $r=u!$ we conclude that there exists $w\in\Nn$ such that $v!=rw$. Hence, in the group $G/M$, we have 
			$(gnM)^{v!}=(gnM)^{rw}=(g^r M)^w=g^{rw} M$ and hence $(gn)^{v!}\in g^{v!} M$. 
		\end{proof}
		
		\begin{thm} \label{thm:basicbound}
			Suppose $\GGG_{A,\ell}$ is connected for all $\ell\in\Ll$ and that $A$ satisfies condition \condH. Let $L\subset\Ll$.
			\begin{enumerate}
				\item[(a)] We have $[\GGG_{A,\ell}^{(q)}:G_{A,\ell}(K)]\le 2c(A/K,L)$ for all $\ell\in L$.
				\item[(b)] We have $y^{2\cdot c(A/K,L)!}\in G_{A,L}(K)$ for all
				$y\in\GGG_{A,L}^{(q)}$.
			\end{enumerate}
		\end{thm}
		
		\begin{proof}
			By our assumptions we have $\GGG_{A,\ell}=\GGG_{A,\ell}^\circ$ for every $\ell\in\Ll$.
			Let $m_{A,L}'$ be the restriction of $m_{A,L}$ to
			$\prod_{\ell\in L}\GGG_{A,\ell}(\Zz_\ell)$.
			Note that by the connectedness assumption
			$\prod_{\ell\in L}\KKK_{A,\ell}(\Zz_\ell)$ is the kernel of $m_{A,L}'$.
			Furthermore
			$\GGG_{A,L}^{(q)}=(m_{A,L}')^{-1}(\langle q_L\rangle)$.
			Thus we have
			$\GGG_{A,L}^{(q)}\supset\prod_{\ell\in L}\KKK_{A,\ell}(\Zz_\ell)$
			and both are normal subgroups of
			$\prod_{\ell\in L}\GGG_{A,\ell}(\Zz_\ell)$.
			We claim that
			\begin{align}\label{eq:claimsquare}
				\GGG_{A,L}^{(q)}
				=G_{A,L}(K)\cdot\prod_{\ell\in L}\KKK_{A,\ell}(\Zz_\ell).
			\end{align}
			Let $y\in\GGG_{A,L}^{(q)}$. Then $m_{A,L}(y)\in\langle q_L\rangle$
			and from $m_{A,L}(G_{A,L}(K))=\langle q_L\rangle$ we conclude
			that there exists $g\in G_{A,L}(K)\subset\GGG_{A,L}^{(q)}$
			such that $m_{A,L}(g)=m_{A,L}(y)$.
			Thus there exists
			$x\in\ker(m_{A,L}')=\prod_{\ell\in L}\KKK_{A,\ell}(\Zz_\ell)$
			such that $y=gx$. This finishes up the proof of the claim.
			
			Let $\ell\in L$. We let
			$S_\ell=G_{A,\ell}(K)\cap\KKK_{A,\ell}(\Zz_\ell)$.
			Equation \eqref{eq:claimsquare}, applied with $L=\{\ell\}$, implies that
			$G_{A,\ell}(K)\cdot\KKK_{A,\ell}(\Zz_\ell)=\GGG_{A,\ell}^{(q)}$, hence
			\begin{align}\label{eq:claimsquare2}
				[\GGG_{A,\ell}^{(q)}:G_{A,\ell}(K)]
				=[\KKK_{A,\ell}(\Zz_\ell):S_\ell].
			\end{align}
			By Lemma \ref{lemm:boundK} we have
			\begin{align}\label{eq:claimsquare-xxx}
				\KKK_{A,\ell}(\Zz_\ell)\triangleright
				\GGG_{A,\ell}^{\geo,\circ}(\Zz_\ell)
				\ \mbox{and}\
				[\KKK_{A,\ell}(\Zz_\ell):
				\GGG_{A,\ell}^{\geo,\circ}(\Zz_\ell)]\le 2.
			\end{align}
			Furthermore
			\begin{align}\label{eq:claimsquare-yyy}
				\GGG_{A,\ell}^{\geo,\circ}(\Zz_\ell)\triangleright
				G_{A,\ell}^\geo(K_A^\ind)
				\ \mbox{and}\
				[\GGG_{A,\ell}^{\geo,\circ}(\Zz_\ell):
				G_{A,\ell}^\geo(K_A^\ind)]\le c(A/K,L)
			\end{align}
			by Lemma \ref{lemma:CHTplusind} and Definition and Remark \ref{defrem:constant}.
			From \eqref{eq:claimsquare-xxx} and \eqref{eq:claimsquare-yyy} it follows that
			\begin{align}\label{eq:claimsquare-zzz}
				[\KKK_{A,\ell}(\Zz_\ell):G_{A,\ell}^\geo(K_A^\ind)]
				\le 2c(A/K,L).
			\end{align}
			Moreover $G_{A,\ell}^\geo(K_A^\ind)\subset S_\ell$ and thus
			\eqref{eq:claimsquare-zzz} and \eqref{eq:claimsquare2} imply Part (a).
			
			To prove Part (b) let $y=(y_\ell)_{\ell\in L}\in\GGG_{A,L}^{(q)}$.
			By \eqref{eq:claimsquare} there exist
			$g=(g_\ell)_{\ell\in L}\in G_{A,L}(K)$ and
			$x=(x_\ell)_{\ell\in L}\in\prod_{\ell\in L}\KKK_{A,\ell}(\Zz_\ell)$
			such that $y=gx$. Let $\ell\in L$.
			
			The groups $\KKK_{A,\ell}(\Zz_\ell)$ and
			$\GGG_{A,\ell}^{\geo,\circ}(\Zz_\ell)$ are normal in
			$\GGG_{A,\ell}(\Zz_\ell)$: the former is the multiplier kernel,
			and the latter is normal because $\GGG_{A,\ell}^{\geo,\circ}$ is the derived group of $\GGG_{A,\ell}$ by Fact \ref{remark:derived}. Moreover, $G_{A,\ell}^\geo(K_A^\ind)$ is characteristic in $\GGG_{A,\ell}^{\geo,\circ}(\Zz_\ell)$ by Lemma \ref{lemma:CHTplusind}~(b), and hence is also normal in $\GGG_{A,\ell}(\Zz_\ell)$.
			
			By \eqref{eq:claimsquare-xxx} and Lemma \ref{lemm:grp:index} we obtain $y_\ell^2=(g_\ell x_\ell)^2\in g_\ell^2\GGG_{A,\ell}^{\geo,\circ}(\Zz_\ell)$. Hence there exists
			$x'_\ell\in\GGG_{A,\ell}^{\geo,\circ}(\Zz_\ell)$
			such that $y_\ell^2=g_\ell^2x'_\ell$. Now \eqref{eq:claimsquare-yyy} and Lemma \ref{lemm:grp:index} give
			$$
			y_\ell^{2\cdot c(A/K,L)!}
			=(g_\ell^2x'_\ell)^{c(A/K,L)!}
			\in g_\ell^{2\cdot c(A/K,L)!}G_{A,\ell}^\geo(K_A^\ind).
			$$
			Thus there exists
			$w\in\prod_{\ell\in L}G_{A,\ell}^\geo(K_A^\ind)$
			such that
			$y^{2\cdot c(A/K,L)!}=g^{2\cdot c(A/K,L)!}w$.
			By Lemma \ref{lemma:CHTplusind}~(a) we see that
			$G_{A,L}^\geo(K_A^\ind)=\prod_{\ell\in L}G_{A,\ell}^\geo(K_A^\ind)$,
			thus $w\in G_{A,L}^\geo(K_A^\ind)$.
			Moreover $G_{A,L}^\geo(K_A^\ind)\subset G_{A,L}(K)$.
			Thus $w\in G_{A,L}(K)$.
			As $g\in G_{A,L}(K)$ it follows that
			$y^{2\cdot c(A/K,L)!}\in G_{A,L}(K)$.
		\end{proof}
		
		\begin{defi} For $L\subset\Ll$ and $s\in\Nn^*$, we define
			$\GGG_{A,L}^{\circ,(q^s)}
			:=\GGG_{A,L}^{(q^s)}\cap
			\prod_{\ell\in L}\GGG_{A,\ell}^\circ(\Zz_\ell)$.
		\end{defi} 
		
		\begin{lemm} \label{lemm:cyclicguy} For $\ell\in \Ll$, $L\subset \Ll$, $j\in\Nn^*$ and $e=[K_A^\conn:K]$ the following hold.
			\begin{enumerate}
				\item[(a)] The index $[\GGG_{A, \ell}^{(q)}:\GGG_{A, \ell}^{\circ, (q^j)}]$ divides $ej$. 
				\item[(b)] For every $g\in \GGG_{A,L}^{(q)}$ we have $g^e\in \GGG_{A, L}^{\circ, (q^e)}$. 
				\item[(c)] If $j$ divides $e$, then  $g^e\in \GGG_{A, L}^{\circ, (q^j)}$. 
			\end{enumerate}
		\end{lemm}
		
		\begin{proof} We consider the group morphism $m_{A,L}: \prod_{\ell\in L} \GGG_{A, \ell}(\Zz_\ell)\to  \prod_{\ell\in L} \Gg_m(\Zz_\ell)$. 
			Now $$ \GGG_{A, L}^{(q^j)}=m_{A,L}^{-1}(\langle q_L^j\rangle)\subset m_{A,L}^{-1}(\langle q_L\rangle)= \GGG_{A, L}^{(q)},$$
			$m_{A,L}(\GGG_{A,L}^{(q^j)})=\langle q_L^j \rangle$ because $G_{A,L}(K)\subset \GGG_{A,L}^{(q)}$, and $m_{A,L}(G_{A,L}(K))=\langle q_L\rangle$. Hence
			$\GGG_{A, L}^{(q)}/\GGG_{A, L}^{(q^j)} \cong \langle q_L \rangle/\langle q_L^j \rangle $ is cyclic group of order dividing $j$. 
			
			The group morphism  $\rho_{A,\ell}$ induces a group isomorphism
			$$\Gal(K_A^\conn/K)\to \GGG_{A, \ell}(\Qq_\ell)/\GGG_{A,\ell}^\circ(\Qq_\ell).$$ Hence $ \GGG_{A, \ell}(\Qq_\ell)/\GGG_{A,\ell}^\circ(\Qq_\ell)$ has 
			cardinality $e$.
			Thus in the chain
			$$\GGG_{A, \ell}^{(q)}\supset \GGG_{A,\ell}^{(q^j)} \supset \GGG_{A, \ell}^{\circ, (q^j)} $$
			the left hand index divides $j$ and the right hand index divides $e$. Part (a) is immediate from that. We now prove (b). 
			Let 
			$g\in \GGG_{A, L}^{(q)}$. Then $g^e\in \prod_{\ell\in L} \GGG_{A,\ell}^\circ(\Zz_\ell)$ and 
			$$m_{A,L}(g^e)=m_{A,L}(g)^e\in \langle q_L^e\rangle.$$
			It follows that $g^e\in  \GGG_{A, L}^{\circ, (q^e)}$. This proves (b). If $j|e$, then $\langle q_L^e\rangle \subset \langle q_L^j\rangle$ and
			$ \GGG_{A, L}^{\circ, (q^e)}\subset  \GGG_{A, L}^{\circ, (q^j)}$. Hence (c) is immediate from (b). 
		\end{proof}
		
		Recall that $k'$ is the exact constant field of $K$. 
		
		\begin{coro} \label{coro:withconstants2} Let $L\subset \Ll$, $e=[K_A^\conn:K]$ and
		$j=[k_A^\conn:k']$.
		Suppose $A/K$ satisfies condition \condH. 
			\begin{enumerate}
				\item[(a)]  We have $[\GGG_{A,\ell}^{(q)}:G_{A, \ell}(K)]\le 2\cdot j\cdot c(A/K, L)$ for all $\ell\in L$.
				\item[(b)]  We have $[\HHH_{A,\ell}^{(q)}:C_{A, \ell}(K)]\le 2\cdot j\cdot c(A/K, L)$ for all $\ell\in L$.
				\item[(c)] We have $g^{2\cdot e\cdot c(A/K, L)!}\in G_{A, L}(K)$ for all $g\in 
				\GGG^{(q)}_{A, L}$. 
				\item[(d)] We have $h^{2\cdot e \cdot c(A/K, L)! }\in C_{A, L}(K)$ for all $h\in 
				\HHH^{(q)}_{A, L}$. 
			\end{enumerate}
		\end{coro}
		
		\begin{proof} Let $\ell\in L$. 
			In the sequel, when applying Theorem \ref{thm:basicbound} over $K_A^\conn$,
			we use the same field $K_A^\ind$ as over $K$.
			This is permissible by Lemma \ref{lemma:CHTplusind},
			since $K_A^\conn\subset K_A^\ind$ and the geometric
			identity components are unchanged under finite extension. Thus, the constant occurring in that application is $c(A/K,L)$.
			Note that $q=\# k'$ and $q^j=\# k_A^\conn$. As $K/k'$ is regular, it follows that $[k_A^\conn K:K]=[k_A^\conn: k']=j$, and that number divides $e$. 
			We have $$\GGG_{A,\ell}^{(q)}\supset \GGG_{A,\ell}^{\circ, (q^j)}\supset G_{A, \ell}(K_A^\conn)
			\ \mbox{and}\ G_{A, \ell}(K_A^\conn)\subset G_{A,\ell}(K).$$ 
			By Theorem \ref{thm:basicbound} applied over $K_A^\conn$, we have  $[\GGG_{A,\ell}^{\circ, (q^j)} : G_{A, \ell}(K_A^\conn)]\le 2c(A/K, L)$. 
			By Lemma \ref{lemm:cyclicguy} (a) $[\GGG_{A,\ell}^{(q)}: \GGG_{A,\ell}^{\circ, (q^j)}]\le ej$.
			Thus
			$[\GGG_{A,\ell}^{(q)}: G_{A,\ell}(K_A^\conn)]\le 2ej c(A/K, L)$ (*). 
			On the other hand 
			$$[\GGG_{A,\ell}^{(q)}:G_{A, \ell}(K_A^\conn)]=[\GGG_{A,\ell}^{(q)}:G_{A, \ell}(K)][G_{A, \ell}(K): G_{A,\ell}(K_A^\conn)]$$
			and $[G_{A, \ell}(K): G_{A,\ell}(K_A^\conn)]=[K_A^\conn: K]=e$. Thus,
			$$[\GGG_{A,\ell}^{(q)}:G_{A, \ell}(K)]=\frac{[\GGG_{A,\ell}^{(q)}:G_{A, \ell}(K_A^\conn)]}{e}\buildrel (*)\over \le 2j c(A/K,L).$$
			This establishes Part (a). 
			
			The group $\HHH_{A,\ell}^{(q)}$ is contained in $\GGG_{A,\ell}^{(q)}$ and $C_{A, \ell}(K)= \HHH_{A,\ell}^{(q)}\cap G_{A, \ell}(K)$. Hence
			$[\HHH_{A,\ell}^{(q)}:C_{A, \ell}(K)]\le [\GGG_{A,\ell}^{(q)}:G_{A, \ell}(K)]$, and thus Part (b) follows from Part (a). 
			
			We prove Part (c). Let $g\in \GGG^{(q)}_{A, L}$. By Lemma \ref{lemm:cyclicguy} (c) $g^e\in \GGG_{A,L}^{\circ, (q^j)}$. 
			Theorem
			\ref{thm:basicbound} applied over $K_A^\conn$ implies that $g^{2e\cdot c(A/K, L)!} =(g^e)^{2\cdot c(A/K, L)!}\in G_{A, L}(K_A^\conn)\subset G_{A,L}(K)$. This proves Part (c). 
			
			Let $h\in \HHH^{(q)}_{A, L}$. Part (c) gives $h^{2\cdot e \cdot c(A/K, L)! }\in  \HHH^{(q)}_{A, L}\cap G_{A,L}(K)=C_{A,L}(K)$. This establishes Part (d). 
		\end{proof}
		
		\begin{rema} The index $j=[k_A^\conn:k']$ divides $[k_A^\conn:k]$. Hence the conclusions of Corollary \ref{coro:withconstants2} hold true with $j$ replaced by
		$[k_A^\conn:k]$. 
		\end{rema} 
		
		\begin{coro} \label{coro:withconstants3} Suppose $A/K$ satisfies condition \condH. 
			Let $\gamma(2g_A)$ be the constant from Lemma \ref{lemm:Kconnbound}. Let $C_{2g_A}$ be the constant from Fact \ref{fact:CHT}. These constants only depend on $g_A$. There exists a finite set $\Sigma\subset \Ll$ (depending on $A/K$) with the following properties. 
			\begin{enumerate}
				\item[(a)]  We have $[\GGG_{A,\ell}^{(q)}:G_{A, \ell}(K)]\le 2\cdot \gamma(2g_A)\cdot C_{2g_A}$ for all $\ell\in \Ll\setminus \Sigma$.
				\item[(b)]  We have $[\HHH_{A,\ell}^{(q)}:C_{A, \ell}(K)]\le 2\cdot \gamma(2g_A) \cdot C_{2g_A}$ for all $\ell\in  \Ll\setminus \Sigma$.
				\item[(c)] We have $g^{2\cdot \gamma(2g_A)\cdot C_{2g_A}!}\in G_{A,  \Ll\setminus \Sigma}(K)$ for all $g\in 
				\GGG^{(q)}_{A,  \Ll\setminus \Sigma}$. 
				\item[(d)] We have $h^{2\cdot \gamma(2g_A) \cdot C_{2g_A}! }\in C_{A,  \Ll\setminus \Sigma}(K)$ for all $h\in 
				\HHH^{(q)}_{A,  \Ll\setminus \Sigma}$. 
			\end{enumerate}
		\end{coro}
		
		\begin{proof}
		Let $\Sigma$ be the finite subset of $\Ll$ from Lemma \ref{lemma:CHTplusind}. We apply 
			Corollary \ref{coro:withconstants2} with $L=\Ll\setminus \Sigma$, and note that $j\mid \gamma(2g_A)$ and $e\mid \gamma(2g_A)$ by Lemma \ref{lemm:Kconnbound} and $C(A/K, L) \le C_{2g_A}$ by Definition and Remark \ref{defrem:constant}. 
		\end{proof}
		
		\section{A new case of the Geyer–Jarden Conjecture}
		
		The exact constant field of $K$ has cardinality $q = p^r$ for an integer $r \geq 1$. For $\ell \in \Ll$, let $p_\ell$ be the image of the prime number $p$ through the unique ring morphism $\Zz \to \Zz_\ell$. If $\alpha_\ell \in \Zz_\ell^\times$, let $\ord_{\Ff_\ell^\times}(\alpha_\ell)$ be the order of its image $\overline{\alpha}_\ell$ in the multiplicative group $\Ff_\ell^\times$ through the reduction map $\Zz_\ell \to \Ff_\ell$. That reduction map yields another $\mathrm{red}_\ell: \GL_{T_\ell(A)}(\Zz_\ell)\to \GL_{T_\ell(A)}(\Ff_\ell)$, and we let $\orho_{A,\ell}$ be $\mathrm{red}_\ell \circ \rho_{A,\ell}$ and $\oG_{A,\ell}(K)$ its image. Let $\oE_{A,\ell}(K)$ be the set of those $\sigma\in \oG_{A,\ell}(K)$ whose action on the $2g_A$-dimensional $\Ff_\ell$-vector space $T_\ell(A)\otimes_{\Zz_\ell}\Ff_\ell=A[\ell](\oK)$ has $1$ among its eigenvalues.
		
		\begin{lemm} \label{lemm:GJ2} If $A/K$ satisfies \condH, let $c_1 := 2 [K_A^\conn:K] c(A/K,\Ll)!$ from Corollary \ref{coro:withconstants2} (d). For all $\ell \in \Ll$, we have $\# \oE_{A,\ell}(K)/\#\oG_{A,\ell}(K) \leq 2g_A/\ord_{\Ff_\ell^\times}(p_\ell^{c_1 r})$.
		\end{lemm}
		\begin{proof} Note that $(q_\ell \Id_{T_\ell(A)})_{\ell \in \Ll} \in \HHH^{(q)}_{A, \Ll}$. Since $A/K$ satisfies \condH, Corollary \ref{coro:withconstants2} (d) implies that for each $\ell \in \Ll$ we have $p_\ell^{c_1r}\Id_{T_\ell(A)}\in C_{A,\ell}(K) \subset G_{A,\ell}(K)$. Let $\oU_\ell$ be the subgroup of $\Ff_\ell^\times$
		generated by $\op_\ell^{c_1 r}$, viewed as a group of scalar matrices. For $u \in \oU_\ell$ and $\sigma \in \oG_{A,\ell}(K)$, we have $u\sigma \in \oE_{A,\ell}(K)$ if and only if $\det(\sigma - u^{-1} \Id_{A[\ell]})=0$. As the characteristic polynomial $\det(\sigma-T\Id_{A[\ell]})\in \Ff_\ell[T]$ has at most $2g_A$ zeros, $\# \left( \oE_{A,\ell}(K)\cap \oU_\ell \sigma\right)\le 2g_A$. Let $m := \#\oG_{A, \ell}(K)/\# \oU_{\ell}$. A system of coset representatives $\{\sigma_1,\cdots, \sigma_m\}$ for the action of $\oU_\ell$ on $\oG_{A,\ell}(K)$ yields a partition $\oG_{A,\ell}(K)=\coprod_{i=1}^m \oU_\ell \sigma_i$. Thus,
			$$
			\# \oE_{A,\ell}(K)\leq \frac{\#\oG_{A, \ell}(K)}{\# \oU_{\ell}}2g_A= \frac{\#\oG_{A, \ell}(K)}{\ord_{\Ff_\ell^\times}(p_\ell^{c_1 r})}2g_A,
			$$
			which concludes.
		\end{proof}
		
		\begin{rema}\label{Remark:Convergent_Series}
			For any integer $j \geq 1$, let $\Ll_n := \{ \ell \in \Ll : \ord_{\Ff_\ell^\times}(\op_\ell^j) = n \}$ and $\Dd_n := \{ \ell \in \Ll : (\op_\ell^j)^n = 1 \}$. Then the sets $(\Ll_n)_{n \geq 1}$ partition $\Ll$ and for each $n \geq 1$ we have $\Ll_n \subset \Dd_n$. Now, $\ell \in \Dd_n$ is equivalent to $\beta \ell = (p^j)^n - 1$ for an integer $\beta \geq 1$ and such $\beta$ is determined by $\ell$. Hence, $\#\Dd_n = \omega((p^j)^n-1)$, the number of distinct prime divisors of $(p^j)^n-1$. As $2^{\omega((p^j)^n-1)} \leq (p^j)^n$, $\omega((p^j)^n-1)/n^3 \leq \log_2(p^j)/n^2$. Hence, for any integer $e \geq 3$, we have
			$$
			\sum_{\ell \in \Ll} \frac{1}{\ord_{\Ff_\ell^\times}(p_\ell^j)^e} \leq \sum_{n=1}^\infty \sum_{\ell \in \Dd_n} \frac{1}{n^3} \leq \log_2(p^j)\sum_{n=1}^\infty \frac{1}{n^2} < \infty.
			$$
		\end{rema}
		
		\begin{thm}\label{gjthm}  Let $K$ be a finitely generated field of positive characteristic, $A/K$ an Abelian variety and $e\ge 3$ an integer. Assume $A$ satisfies condition \condH. For all $\sigma\in \Gal(K)^e$ outside a subset of $\Gal(K)^e$ of measure zero there are only finitely many  prime numbers $\ell$ with $A(\oK(\sigma))[\ell]\neq 0$.
		\end{thm}
		
		\begin{proof} It suffices to consider $\ell \in \Ll$ and for each of those, let $S_{e,\ell}$ be the set of $\sigma \in \Gal(K)^e$ such that $A(K^\sep(\sigma))[\ell] \neq 0$. Let
			$$
			S := \bigcap_{n=1}^\infty \bigcup_{\substack{\ell \in \Ll \\ \ell \geq n}} S_{e,\ell} = \{ \sigma \in \Gal(K)^e : \sigma \text{ belongs to infinitely many } S_{e,\ell} \}.
			$$
			For each $\ell \in \Ll$, $S_{e,\ell} \subset \orho_{A,\ell}^{-1}(\oE_{A,\ell}(K))^e$ and by Lemma \ref{lemm:GJ2} we have
			\begin{equation}\label{Bound_on_measure}
				\mu(S_{e,\ell}) \leq \mu(\oE_{A,\ell})^e = \left( \frac{\# \oE_{A,\ell}}{\# \oG_{A, \ell}(K)} \right)^e \leq \frac{(2g_A)^e}{\ord_{\Ff_\ell^\times}(p_\ell^{c_1 r})^e}.
			\end{equation}
			
			Equation
			\eqref{Bound_on_measure} followed by Remark \ref{Remark:Convergent_Series} yields
			$$
			\sum_{\ell \in \Ll} \mu(S_{e,\ell}) \leq (2g_A)^e\sum_{\ell \in \Ll} \frac{1}{\ord_{\Ff_\ell^\times}(p_\ell^{c_1 r})^e} < \infty.
			$$
			By the first Borel-Cantelli Lemma \cite{Fried_Jarden_2008}*{Lemma 18.3.5 (a)}, $\mu(S) = 0$. Hence, outside a $\mu$-measure zero subset of $\Gal(K)^e$ there are only finitely many $\ell \in \Ll$ with $A(K^\sep(\sigma))[\ell] \neq 0$. We conclude since for all $\ell \in \Ll$, we have $A(K^\sep(\sigma))[\ell] = A(\oK(\sigma))[\ell]$.
		\end{proof}
		
		\subsection*{Tool and computational resource disclosure}
			The authors used AI assisted tools only for proofreading parts of the manuscript. These tools were not used to generate the mathematical content, develop proofs or to write the manuscript. The authors take full responsibility for the content of the paper.

\input{bib}
	\end{document}

%% file: math.tex
\newtheorem{leer}{}[section]
\newtheorem{thm}[leer]{Theorem}
\newtheorem{conj}[leer]{Conjecture}
\newtheorem{rema}[leer]{Remark} 
\newtheorem{prop}[leer]{Proposition}
\newtheorem{lemm}[leer]{Lemma}

\newtheorem{coro}[leer]{Corollary}

\newtheorem{defi}[leer]{Definition}
\newtheorem{defirema}[leer]{Definition and Remark}

\newtheorem{fact}[leer]{Fact}

\newtheorem{intthm}{Theorem}

\newtheorem{intcoro}[intthm]{Corollary}

\newtheorem{intprop}[intthm]{Proposition}

\newcommand{\ol}[1]{\overline{#1}}

\newcommand{\plim}{\mathop{\varprojlim}\limits}

\newcommand{\mf}{\mathfrak}

\newcommand{\Spec}{\mathrm{Spec}}
\newcommand{\geo}{\mathrm{geo}}
\newcommand{\ind}{\mathrm{ind}}
\newcommand{\Id}{\mathrm{Id}}
\newcommand{\End}{\mathrm{End}}
\newcommand{\Lie}{\mathrm{Lie}}
\newcommand{\alg}{\mathrm{alg}}
\newcommand{\GSp}{\mathrm{GSp}}
\newcommand{\diag}{\mathrm{diag}}

\newcommand{\LL}{{\mathcal L}}

\newcommand{\AAA}{\mathscr{A}}

\newcommand{\DDD}{\mathscr{D}}

\newcommand{\GGG}{\mathscr{G}}
\newcommand{\HHH}{\mathscr{H}}

\newcommand{\KKK}{\mathscr{K}}
\newcommand{\LLL}{\mathscr{L}}
\newcommand{\MMM}{\mathscr{M}}

\newcommand{\SSS}{\mathscr{S}}
\newcommand{\TTT}{\mathscr{T}}

\newcommand{\ZZZ}{\mathscr{Z}}

\newcommand{\Aa}{{\mathbb{A}}}

\newcommand{\Cc}{{\mathbb{C}}}
\newcommand{\Dd}{{\mathbb{D}}}

\newcommand{\Ff}{{\mathbb{F}}}
\newcommand{\Gg}{{\mathbb{G}}}

\newcommand{\Ll}{{\mathbb{L}}}

\newcommand{\Nn}{{\mathbb{N}}}

\newcommand{\Qq}{{\mathbb{Q}}}

\newcommand{\Zz}{{\mathbb{Z}}}

\newcommand{\mfc}{{\mf{c}}}

\newcommand{\mfg}{{\mf{g}}}

\newcommand{\mfs}{{\mf{s}}}

\newcommand{\mfz}{{\mf{z}}}

\newcommand{\mfG}{{\mf{G}}}

\newcommand{\ok}{{\ol{k}}}

\newcommand{\op}{{\ol{p}}}
\newcommand{\oq}{{\ol{q}}}

\newcommand{\oE}{{\ol{E}}}

\newcommand{\oG}{{\ol{G}}}

\newcommand{\oI}{{\ol{I}}}

\newcommand{\oK}{{\ol{K}}}

\newcommand{\oU}{{\ol{U}}}

\newcommand{\GL}{\mathrm{GL}}

\newcommand{\Gal}{\mathrm{Gal}}
\newcommand{\rk}{\mathrm{rk}}

\newcommand{\im}{\mathrm{im}}
\newcommand{\chara}{\mathrm{char}}

\newcommand{\ord}{{\mathrm{ord}}}